\documentclass[3p,letterpaper]{elsarticle}
\pdfoutput=1

\usepackage{times} 
\usepackage{amsmath,amssymb,amsfonts}
\usepackage[latin1]{inputenc}
 \usepackage[english]{babel}
\usepackage[cyr]{aeguill}
\usepackage{xcolor} 
\usepackage{geometry} 
\usepackage[pdftex,a4paper,hyperindex=true,hyperfigures=true,bookmarks=true,pdftitle={TITLE},pdfsubject={},pdfauthor={AUTHOR},pdfkeywords={KEYWORDS},pdfpagemode=UseOutlines,colorlinks=true,urlcolor=blue,citecolor=red,linkcolor=blue]{hyperref}
 \ifx\pdfoutput\undefined 
    \usepackage[dvips]{graphicx}
\else                    
 \usepackage[]{graphicx}
\fi
\usepackage{float} 
\DeclareGraphicsExtensions{.jpg,.mps,.pdf,.png} 
\usepackage[sf]{caption2}
\usepackage[figuresright]{rotating}

\newcommand{\review}[1]{{#1}}

\makeindex  
\makeglossary  

\newcommand{\RR}{{\mathbb R}}

\newcommand{\ep}{{\varepsilon}}
\newcommand{\real}{{\mathbb R}}
\newcommand{\R}{{\mathbb R}}

\newcommand{\ZZ}{{\mathbb Z}}

\newcommand{\EX}[1]{\mathbb{E}\!\left[#1\right]}
\newcommand{\la}{\lambda}

\newcommand{\al}{\alpha}

\newcommand{\Hmin}{{h^{\min}}}
\newcommand{\pL}{\ell^{(p)}}
\newcommand{\cpL}{\bar\ell^{(p)}}
\newcommand{\TDFA}{{T_{d}}}
\newcommand{\TMFD}{{T}_{mfd}}
\newcommand{\LMFD}{\mathcal{L}_{mfd}}

\def\w{\omega_0}
\def\ww{\omega_1}
\def\cjk{c_{j,k}}

\newcommand{\BE}{\begin{equation}}
\newcommand{\EE}{\end{equation}}

\newtheorem{lem}{Lemma}
\newtheorem{coro}{Corollary}
\newtheorem{Theo}{Theorem}
\newtheorem{prop}{Proposition}
\newtheorem{defi}{Definition}
\newcommand{\BP}{\begin{prop}}
\newcommand{\EP}{\end{prop}}
\newcommand{\BC}{\begin{coro}}
\newcommand{\EC}{\end{coro}}
\newcommand{\BL}{\begin{lem}}
\newcommand{\EL}{\end{lem}}
\newcommand{\BD}{\begin{defi}}
\newcommand{\ED}{\end{defi}}
\newcommand{\BT}{\begin{Theo}}
\newcommand{\ET}{\end{Theo}}

\begin{document}
\begin{frontmatter}

\title{\bf\boldmath$p$-exponent and $p$-Leaders, Part II:\\ Multifractal Analysis.  \\
Relations to Detrended Fluctuation Analysis.}

\author[lyon]{R. Leonarduzzi}
\ead{roberto.leonarduzzi@ens-lyon.fr}

\author[toul]{H. Wendt}
\ead{herwig.wendt@irit.fr}

\author[lyon]{P. Abry\corref{mycorraut}}
\ead{patrice.abry@ens-lyon.fr}

\author[cret]{S. Jaffard}
\ead{jaffard@u-pec.fr}

\author[mars]{C. Melot}
\ead{melot@cmi.univ-mrs.fr}

\author[lyon]{S. G. Roux}
\ead{stephane.roux@ens-lyon.fr}

\author[arge]{M.~E. Torres}
\ead{metorres@santafe-conicet.gov.ar}

\address[lyon]{%
Signal, Systems and Physics, Physics Dept., CNRS UMR 5672, Ecole Normale Sup\'erieure de Lyon, Lyon, France}
\address[toul]{%
IRIT, CNRS UMR 5505, University of Toulouse,  France}
\address[cret]{%
Universit\'e Paris Est, Laboratoire d'Analyse et de Math\'ematiques Appliqu\'ees, CNRS UMR 8050, UPEC,  Cr\'eteil, France}
\address[mars]{%
Aix Marseille Université, CNRS, Centrale Marseille, I2M, UMR 7373, 13453 Marseille, France}  
\address[arge]{%
Consejo Nacional de Investigaciones Cient\'{\i}ficas y T\'ecnicas, Universidad Nacional de Entre Ríos, Argentina}

\cortext[mycorraut]{Corresponding author. Tel: +33 47272 8493. Postal address: CNRS, Laboratoire de Physique Ecole Normale Superieure de Lyon 46, allée d'Italie, F-69364, Lyon cedex 7, France}

\begin{abstract} 
Multifractal analysis studies signals, functions, images or fields via the fluctuations of their local regularity along time or space, which capture crucial features of their temporal/spatial dynamics.
It has become a standard signal and image processing tool and is commonly used in numerous applications of different natures. 
In its common formulation, it relies on the H\"older exponent as a measure of local regularity, which is by nature restricted to positive values and can hence be used for locally bounded functions only. 
In this contribution, it is proposed to replace the H\"older exponent with a collection of novel exponents for measuring local regularity, the $p$-exponents. One of the major virtues of $p$-exponents is that they can potentially take negative values. 
The corresponding wavelet-based multiscale quantities, the $p$-leaders, are constructed and shown to permit the definition of a new multifractal formalism, yielding an accurate practical estimation of the multifractal properties of real-world data. 
Moreover, theoretical and practical connections to and comparisons against another multifractal formalism, referred to as multifractal detrended fluctuation analysis, are achieved. 
The performance of the proposed $p$-leader multifractal formalism is studied and compared to previous formalisms using synthetic multifractal signals and images, illustrating its theoretical and practical benefits. 
The present contribution is complemented by a companion article studying in depth the theoretical properties of $p$-exponents and the rich classification of local singularities it permits.
\end{abstract}

\end{frontmatter}

\section{Introduction}

\noindent {\bf Context: Scale invariance and multifractal analysis.} \quad 
The paradigm of scale-free dynamics, scale invariance, or scaling, has been shown to be relevant to model empirical data from numerous real-world applications of very different natures, ranging from Neurosciences \cite{Goldberger2002,Werner10,He14,Ciuciu12}, heart rate variability \cite{Ivanov2007,Kiyono2006,Doret2011}, bones textures \cite{benhamou2001bone} to physics (turbulence \cite{m74}), geophysics (rainfalls \cite{Foufoula-Georgiou1994}, \review{wind \cite{telesca2011analysis}, earthquakes \cite{telesca2006}}), finance \cite{Mandelbrot1997,Mandelbrot1999,Lux2007}, Internet traffic \cite{Abry2002}, art investigations \cite{AJWVanGogh2013}, \review{music \cite{telesca2011revealing}} \ldots
In essence, scale invariance is associated to signals whose temporal dynamics involve a wide range (continuum) of time/space scales, rather than being dominated by one or a few specific time/space scales.

Multifractal analysis provides a generic framework for studying scaling properties both theoretically and practically.
Multifractal analysis describes the dynamics of the fluctuations along time or space of the local regularity, $h(x_0)$, of the signal (or image) $X$ around location $x_0$. 
However, though theoretically grounded on local or pointwise quantities, multifractal analysis actually does not aim to produce estimates of $h$ for all $x$. Instead, it  provides practitioners with a global, geometrical and statistical description of the fluctuations of $h$ across all time or space locations. 
This \emph{global} information is encoded in the so-called \emph{multifractal spectrum} ${\cal D}(h) $ (also termed singularity spectrum) which is defined as the fractal dimension of the set of points $x$ where the pointwise regularity exponent $h (x)$ takes the value $h$.

For real-world data, the estimation of the multifractal spectrum does not rely directly on local estimates of $h$ but instead involves a thermodynamics-inspired formalism, the so-called \emph{multifractal formalism}.
A multifractal formalism first requires the choice of \emph{multiresolution quantities}, hereafter labelled $T_X(a,k)$, i.e., quantities that depend jointly on the position $k$ and the analysis scale $a$. 
Scale invariance is defined as the power-law behavior, with respect to the analysis scale $a$, of the time average of the $q-$th power of these $T_X(a,k)$:
\begin{equation}
\label{equ-si}
S(a,q) = {1 \over n_a} \sum_{k=1}^{n_a}  |T_X(a,k)|^q  \simeq a^{\zeta(q)}.
\end{equation}
The quantity $S(a,q) $ is referred to as the {\em structure function} associated with the quantities $T_X(a,k)$, and the function $\zeta(q)$ is termed the \emph{scaling function}. 
The scaling function $\zeta(q)$ can be related to functional space signal or image characterization, cf. e.g., \cite{Jaffard2004,Bergou}.
Most importantly, the multifractal formalism relates the  Legendre transform  ${\cal L}  (h)$ of $\zeta(q) $ to ${\cal D}(h)$ via the formula
 practical estimation of ${\cal D}(h)$ can be obtained via the: 
\begin{equation}
\label{equ-tl}
{\cal D}(h) \leq {\cal L} (h)    : =  \inf_{q \in \real} (d+qh - \zeta(q))
\end{equation}
(cf. e.g., \cite{Jaffard2004,Abel}), and strongly ties the paradigm of \emph{scale invariance} to multifractal properties and the spectrum ${\cal D}(h)$ and  provides a way of estimating ${\cal D}(h)$ in practice. 

For a number of synthetic processes with known ${\cal D}(h)$, it can be shown that the inequality \eqref{equ-tl} turns out to be an equality.
It is also well-known that earlier formulations of multifractal formalisms relying on increments or wavelet coefficients as multiresolution quantities yield poor or incorrect estimates of the multifractal spectrum. 
Commonly used formalisms are relying on Wavelet Transform Modulus Maxima (WTMM) \cite{muzy1993multifractal} or MultiFractal Detrended Fluctuation Analysis (MFDFA) (cf. \cite{kantelhardt2002multifractal} for the founding article, \review{and \cite{gu2010detrending,schumann2011multifractal} for related developments}). 
Recently, it has been shown that a formalism based on Wavelet Leaders provides a theoretically well grounded and practically robust and efficient framework for the estimation of ${\cal D}(h)$ \cite{Wendt:2007:E,WENDT:2009:C,MandMemor,Abel}.
\vskip4mm

\noindent{\bf Multiresolution quantities and regularity exponents.\quad}
An often overlooked issue lies in the fact that the choice of multiresolution quantity $T_X$ is crucial, both theoretically and practically: 
The multifractal formalism \eqref{equ-tl} relies on the assumption that the pointwise exponent $h_X$ can be theoretically recovered through  local log-log regressions of the multiresolution quantities
\begin{equation}
\label{loglog}
h(x) = \liminf_{a \rightarrow 0} \frac{\log (T_X (a, k(x))) }{\log a }
\end{equation}
(where $T_X (a, k(x))$ denotes the value taken by $T_X (a, k)$ at the location $k(x)$ closest to $x$).
Therefore, choosing the multiresolution quantity $T_X (a, k)$ on which the formalism is based also amounts to choosing (and thus changing) the definition of local regularity.  
The predominantly if not exclusively employed notion of local regularity $h$ is the H\"older exponent (cf. Section~\ref{sec:holex}). 
Typical examples for $T_X$ that have been associated with H\"older exponents are increments, oscillations, continuous and discrete wavelet coefficients, wavelet transform modulus maxima (WTMM) \cite{muzy1993multifractal,muzyetal94,ard02,Arneodo2003a}, or wavelet leaders \cite{Jaffard2004,JAFFARD:2006:A}. 
Note, however, that oscillations and wavelet leaders are the only multiresolution quantities for which it has been shown theoretically that \eqref{loglog} holds and that it actually characterizes the H\"older exponent \cite{Jaffard2004}. 
It has been shown that increments and wavelet coefficients are associated to another (and weaker) notion of regularity \cite{JAFFARD:2010:A,JAFFARD:2006:A,MandMemor}, thus partially explaining their poor estimation performance. 
So far, there are no theoretical results available connecting WWTM or MFDFA to H\"older regularity. 
It is shown here in Section~\ref{sec-dfa} below that the latter, MFDFA, is not associated to H\"older regularity but to a different local regularity measure.
At present, wavelet leaders constitute the core of the current state-of-the-art multifractal formalism associated to H\"older regularity \cite{Wendt:2007:E,WENDT:2009:C}. 

However, the H\"older exponent based measure of regularity suffers from the fundamental drawback that it is, by definition, positive and can hence not characterize negative regularity. Multifractal analysis is therefore restricted to functions whose local regularity is positive everywhere. 
This represents a severe limitation since data are often found to contain negative regularity in applications (cf. e.g., \cite{MandMemor,Bergou}).
\vskip4mm

\noindent{\bf Goals, contributions and outline.\quad}
In this context, the present contribution aims to overcome this restriction for the application of multifractal analysis by making use of new pointwise regularity exponents, the $p$-exponents, and the corresponding multiresolution quantities, the $p$-leaders, which have been proposed in the mathematical literature recently \cite{JaffMel}.
The $p$-exponents offer a more versatile framework than H\"older exponents and notably enable to theoretically define and measure \emph{negative regularity}.
In the companion paper \cite{PART1}, $p$-exponent regularity and $p$-leaders have been theoretically defined, studied and characterized. 
The first goal of the present contribution is the development of the corresponding multifractal formalism. 
The second goal is to compare the $p$-leader multifractal formalism theoretically and practically against wavelet leaders and against MFDFA. 

After having formulated a wavelet-based definition of $p$-exponents and $p$-leaders and having briefly restated their key properties (Section \ref{sec:pExp}), the $p$-leader multifractal formalism is derived in Section \ref{sec-mfa}, and explicit estimation procedures for the multifractal spectrum that take into account the discrete and finite-size nature of data are developed.
Section~\ref{sec-dfa} starts with a brief description of the original formulation of DFA and recalls its multifractal extension, the MFDFA method.
It establishes, for the first time, clear theoretical and practical connections between MFDFA and $p$-exponents and points out the theoretical and practical consequences of the conceptual differences between MFDFA and $p$-leaders.
Section~\ref{sec-num} provides a numerical study and validation of the proposed $p$-leader multifractal formalism and compares it with the wavelet leader and MFDFA formalisms in terms of estimation performance, illustrating the clear practical benefits of the proposed novel multifractal formalism. 

A {\sc Matlab} toolbox implementing estimation procedures for performing $p$-multifractal analysis will be made publicly available at the time of publication.

\section{\boldmath Wavelet based definition of pointwise $p$-exponents}
\label{sec:pExp}  

\subsection{Wavelet characterization of uniform regularity}

\subsubsection{Discrete wavelet decomposition}

Wavelet analysis has already been massively used for the study of the regularity of a function  \cite{Jaf91, muzyetal94,jmf2,Riedi03,Jaffard2004}.
Let $\{\psi^{(i)}(x)\}_{ i = 1, \cdots, 2^d-1} $ denote a family of {\em mother wavelets}, consisting of oscillating functions with fast  decay and good joint time-frequency localization. 
Mother wavelets are defined so that the collection of dilated (to scales $2^{j}$) and translated (to space/time location $2^{j}k$) templates
\BE
\label{wavbasbad}  
\Big\{
2^{-dj/2}\psi^{(i)}(2^{-j}x-k), \hspace{5mm} \mbox{ for} \hspace{3mm} i = 1, \cdots 2^d -1, \hspace{3mm}   j\in \ZZ \hspace{3mm} \mbox{ and} \hspace{3mm}
k\in\ZZ^d
\Big\},
\EE
forms an orthonormal basis of $L^2 (\RR^d)$  \cite{Mey90I}.
Without loss of generality, a $d$-dimensional orthonormal wavelet basis is defined here by tensor product of the univariate orthonormal wavelet basis.
Mother wavelets $\psi^{(i)}(x)$ are moreover defined to guarantee a {\em number of vanishing moments}, a strictly positive integer $N_\psi \geq 1$ such that
$
  \int_\R P(x)  \psi^{(i)}(x) dx =  0
$
for any polynomial $P$ of degree strictly smaller than $N_\psi $, while this integral does not vanish for some polynomials of degree $N_\psi$.  
The coefficients  of the $d-$dimensional discrete wavelet transform  of $X$ are defined as
\begin{equation}
\label{eq-wc}
c^{(i)}_{j,k} =    \displaystyle\int_{\R^d} X(x)  \;  2^{-dj}\psi^{(i)}(2^{-j}x-k) \, dx
\end{equation}
where the $L^1$ normalization for the wavelet coefficients is better suited in the context of local regularity analysis, see \cite{Wendt:2007:E,Abel}.

\subsubsection{Wavelet structure function  and scaling function} 

Let $q> 0$. The wavelet structure function is defined as
\begin{equation}
\label{equ-WSF}
S_c(j,q) = 2^{dj} \displaystyle\sum_{  k }\sum_{i=1}^{2^d-1}  \big| c^{(i)}_{j,k} \big|^q, 
\end{equation}
and  the wavelet scaling function as
\BE 
\label{defscalond}  
\forall q >0, \hspace{6mm} 
 \eta(q) =   \displaystyle \liminf_{j \rightarrow - \infty} \;\; \frac{\log \left( S_c(j,q)  \right) }{\log (2^{j})}. \EE 
The wavelet scaling function $\eta(q)$ does not depend on the (sufficiently smooth) wavelet basis which is used, see \cite{JAFFARD:2010:A,Jaffard2004} and Section ~3.1 of companion article \cite{PART1} . 

The wavelet scaling function is of central importance for validating \emph{uniform regularity assumptions} on the function $X$ in  \eqref{equ:p0} below, see the companion paper \cite{PART1} and \cite{Bergou,JAFFARD:2010:A,Abel,MandMemor} for details.
The wavelet scaling function also serves as an important auxiliary quantity in the construction of estimators for finite resolution data, cf. Section \ref{sec:finite}.

\subsection{$p$-leader based definitions of $p$-exponents}  
 \label{usageeta}
 
 \subsubsection{ $p$-leaders}
 
Let us \review{now} define the multiresolution quantities at the heart of the present contribution, the \emph{$p$-leaders} (cf. the companion paper \cite{PART1} for details).
Let $ k = (k_1, \ldots, k_d)$, $j\in\ZZ$ and let the dyadic cubes be indexed by
$$ 
\lambda\; (= \lambda_{j,k}) \; : = \left[ 2^j k_1,  2^j (k_1+1)  \right) \times \cdots \times  \left[ 2^j k_d,  2^j (k_d+1) \right)
$$
and, accordingly: $c^{(i)}_{\lambda} = c^{(i)}_{ j,k}$.  
Furthermore, let $3 \lambda$ denote the cube with the same center as that of the cube $\la$ but $3$ times wider.
\BD 
\label{def:pleaders}
Let $p>0$ . 
The \emph{$p$-leaders} are defined as  \cite{JaffToul, Bergou,JaffMel,PART1}
\begin{equation} 
\label{pleaders} 
 \pL(j,k)  \equiv \pL_\la= \left( \sum_{j' \leq j, \;  \lambda' \subset 3 \lambda} \sum_{i=1}^{2^d-1} \big| c^{(i)}_{\la'} \big|^p \, \, 2^{-d(j-j')} \right)^\frac{1}{p},\quad p>0
\end{equation}
where $j'\leq j$ is the scale associated with the sub-cube $\la'$ of width $2^{j'}$ included in  $3 \lambda$.
\ED
This definition is illustrated in \cite[Fig. 1]{PART1}. 
The sum in \eqref{pleaders} is finite if the condition 
$\eta(p) >0$ is satisfied, see \eqref{equ:p0} below.

 \subsubsection{$p$-exponents}

When $p$-leaders are well defined ($\eta(p)>0$), the  \emph{pointwise $p$-exponent} $h_p(x_0)$ can be defined as follows \cite{JaffToul, Bergou,JaffMel}.
\BD 
\label{def:pexpo}
Let $p>0$ and $\eta(p)>0$. The pointwise $p$-exponent of $X$ at $x_0$ is 
\begin{equation} 
\label{carachqf}  
h_p(x_0)  =  \liminf_{j \rightarrow - \infty}  \frac{ \log  \left(   \pL_{\lambda_j (x_0)} \right)  }{\log (2^{j})}, 
\end{equation}
where $ \la_j (x_0)$ denotes the dyadic cube of width $2^{j} $ that contains $x_0$.
\ED
This definition does not depend on the wavelet as long as the wavelet is smooth enough (cf. \cite{PART1} Section~3.1 and\cite[Chapter 3]{Mey90I}).

In practice, (\ref{carachqf}) essentially means  that
\begin{equation}
\label{carachqfsim}  
\pL_{\lambda_j (x_0)} \sim C 2^{jh_p(x_0)} \, \mbox{when} \; a = 2^{j} \rightarrow 0, 
\end{equation}
which meets the natural prerequisite (cf. \eqref{loglog}) to construct a multifractal formalism. 
This will be the subject of Section \ref{sec-mfa}.

When $p\geq 1$, Definition \ref{def:pexpo} is equivalent to the original definition of $p$-exponents supplied by the $T^p_\alpha (x_0)$  regularity of A. Calder\'on and A. Zygmund \cite{CalZyg}, 
which essentially states that
\begin{equation}  
\label{equ:tpa}
T^{(p)}(a,x_0)= \left( \frac{1}{a^d} \int_{B(0,a)} | X(u+x_0) -P_{x_0}(u+x_0)|^p du  \right)^\frac{1}{p} \sim C a^{h_p(x_0)}
\end{equation}
where $B(x_0, a)$ denotes the ball of radius $r$ centered at $x_0$, $P_{x_0}$ is a polynomial  of degree less than $h_p(x_0)$ and $C,R$ are positive constants, see \cite{JaffMel,JaffToul,Bergou} and \cite[Definition 1]{PART1} for details.
 Note that the advantage of Definition \ref{def:pexpo} is that it is valid for all $p>0$.
Furthermore, note that for the case $p=+\infty$, we recover from \eqref{pleaders} the definition of classical \emph{wavelet leaders}, see, e.g., \cite{Jaffard2004,JAFFARD:2006:A,Wendt:2007:E,WENDT:2009:C}
$$
\ell_\lambda = \ell_\lambda^{(+\infty)} = \sup_{i \in [ 1,\ldots, 2^d-1], \; \lambda ' \subset 3 \lambda} | c^{(i)}_{\lambda'} |
$$
and from \eqref{equ:tpa} the definition of the classical \emph{H\"older exponent} \cite{Jaffard2004}
\begin{equation}
\label{equ:tpe}
|X(a+x_0)-P_{x_0}(a+x_0)| \sim C |a|^{h(x_0)}, \, \, |a| \rightarrow 0.
\end{equation}
\vskip2mm

When $0<p<1$, this wavelet-based Definition  \ref{def:pexpo} of $p$-exponents has been shown to be well-grounded, relevant and useful (cf. \cite{PART1}).

\subsection{ Properties of $p$-exponents}
\label{sec:pproperties}
 The following theoretical key properties  are satisfied by $p$-exponents (details and proofs are given in the companion paper \cite{PART1}).

\noindent{\bf Domain of definition.\quad} $p$-exponents are defined for $p\in(0,p_0]$, where the critical Lebesgue index $p_0$ is defined as%
\footnote{%
Note that \eqref{equ:p0} actually provides a lower bound for the  quantity $p_0 (x_0) $ whose definition is given in \cite[Theorem 1]{PART1} and will be  used here as an operational surrogate. 
}%
\begin{equation}
\label{equ:p0}
p_0 = \sup \{ p:  \eta(p)>0\}.
\end{equation}
In contradistinction,  H\"older exponents are only defined for locally bounded functions, i.e., functions for which $p_0=+\infty$.

\noindent{\bf Negative regularity.\quad} $p$-exponents can take values  down to $-d/p$,  see \cite[Theorem 1]{PART1}. Therefore, $p$-exponents enable the use of { \em negative regularity exponents}. For instance, they enable to model local behaviors of the form $|X(x)|\sim 1/|x-x_0|^\alpha$ for $\alpha<d/p$.

\noindent{\bf\boldmath $p$-exponents for different values of \boldmath$p$.\quad} $p$-exponents for different values of $p$ do not in general coincide but satisfy $ h_p(x_0) \geq h_{p'}(x_0)$ for $p \leq p' \leq p_0$, see \cite[Theorem 1]{PART1}. Therefore, $p$-exponents for different values of $p$ can provide complementary information for the characterization of singularities.

\noindent{\bf Classification of singularities.\quad} $p$-exponents enable a generic classification of singularities, see \cite[Section 4]{PART1} for details and examples.
Let $X$ have a singularity at $x_0$ that has $p$-exponent $h_p(x_0)$. The singularity of $X$ at $x_0$ is \emph{$p$-invariant} if the function $p\to h_p(x_0)$ does not depend on $p$. Moreover,
\begin{itemize}
\item[-] $X$ has a \emph{canonical singularity} at $x_0$ if the $p$-exponent of its fractional integral $X^{(-t)}$ of order $t$ equals  $h_p(x_0)+t$ (cf., \cite[Definition 3-4]{PART1}). 
 Examples are provided by deterministic self-similar functions, e.g., the cusps  $|x-x_0|^\alpha$, $\alpha>0$, see  \cite[Proposition 1]{PART1}. 
An important property of canonical singularities is that they form a subclass of the  $p$-invariant singularities  \cite[Theorem 2]{PART1}.
\item[-] A singularity that is not canonical is called an \emph{oscillating singularity}:
\begin{itemize}
\item[-] An oscillating singularity that is $p$-invariant is termed \emph{balanced}, see  \cite[Definition 7]{PART1}. An example is provided by the ``chirp'' type singularities $|x-x_0|^\alpha\sin(1/|x-x_0|^\beta)$, $\alpha,\beta>0$.
\item[-] An oscillating singularity for which the function $p\to h_p(x_0)$ \emph{does} depend on $p$ is termed \emph{lacunary}. An example is given by the lacunary comb $F_{\alpha,\gamma}(x)$ in  \cite[Eq. (12)]{PART1} and by the random processes in Section \ref{sec:rws}.
\end{itemize}
\end{itemize}

\noindent{\bf H\"older exponent.\quad} The H\"older exponent is given by the $p$-exponent with $p =+\infty$: $h(x_0) \equiv h_{+\infty}(x_0)$.  
However, it is central to realize that $p$-exponents for finite values of $p$ measure a regularity in $X$ that does \emph{not} in general coincide with that captured by the H\"older exponent, cf., \cite[Theorem 1]{PART1}.

When $p\rightarrow\infty$, Condition $\eta(p)>0$ in \eqref{equ:p0} is the counterpart to the condition met in wavelet leader multifractal analysis (cf. \cite{PART1,MandMemor,Bergou}):
\begin{equation}
\label{caracbeswav3hol}  
\Hmin =\liminf_{j \rightarrow - \infty} \;\;  
 \;\;  \frac{ \log 
\left(  \displaystyle  \sup_{ i,k}  \big| c^{(i)}_{ j,k} \big|     \right) }{\log (2^{j})}
> 0. 
\end{equation}
 
\section{Multifractal Analysis}
\label{sec-mfa}

In practice, the values taken by a  pointwise regularity exponent $h(x)$ cannot be extracted point by point  from data $X(x)$. Instead, multifractal analysis  aims  to provide a global description of the fluctuations along time or space of $h(x)$, termed the singularity or \emph{multifractal spectrum}, which can be estimated  from data $X(x)$ by recourse to a so-called \emph{multifractal formalism}. 
In the case of H\"older regularity, such a formula is given by the {wavelet leader multifractal formalism}, cf., e.g., \cite{Jaffard2004,JAFFARD:2006:A,Wendt:2007:E,WENDT:2009:C,Abel}.
The goal of this section is to define and study a multifractal formalism based on $p$-exponents and $p$-leaders and to derive explicit expressions for estimators that can be applied to data in applications.

The $p$-leaders and $p$-exponents have not yet been used for practical multifractal analysis intended for real-world data analysis except in our preliminary works
\cite{Leonarduzzi2014icassp,Leonarduzzi2014CLAIB,Wendt2014EMBC}.

\subsection{ Multifractal $p$-spectrum}
\label{sec:holex}

The \emph{multifractal $p$-spectrum} of $X$ is defined as follows.
\BD The  {multifractal $p$-spectrum} ${\cal D}^{(p)}(h) $ is  the Hausdorff dimension $\dim_H$ of the set of points where the  $p$-exponent takes the value $h$
\begin{equation}
\label{eq:Dh}
{\cal D}^{(p)}(h) = \dim_H \left(  \{ x \in \real^d, h_p(x) = h \}  \right).
\end{equation}
The support of the spectrum   is  the image of the mapping $x\rightarrow h_p(x)$, i.e. the collection of values   of $h$ such that $ \{ x \in \real^d:  h_p(x) = h \} \neq \emptyset $. 
\ED
By convention $\dim (\emptyset ) = -\infty$. 

Because $p$-exponents are necessarily larger than $-d/p$, the  support of ${\cal D}^{(p)}$ is included in  $[  -d/p, +\infty ] $. 
In addition, ${\cal D}^{(p)}$ is further constrained by the following sharper condition:
\BP \label{majpspec} Let $p > 0$ and let $X$ be a function for which $\eta(p)>0$. Then 
\BE 
\label{bound:dh}\forall h \leq 0, \qquad {\cal D}^{(p)}(h) \leq d+hp.  \EE
\EP
The proof is given in \ref{proofs}.  

This condition implies that the multifractal $p$-spectrum is necessarily below a straight line that 
connects points $(-d/p, 0)$ and $(0,d)$.
This is illustrated in Figs. \ref{fig:illMRW} and \ref{fig:CMC2D}. 

For further details on multifractal analysis, notably for the definition of the Hausdorff dimension,
interested readers are referred to \cite{Fal93,jmf2,Riedi03,Jaffard2004}.

\subsection{ $p$-leader multifractal formalism}  
 \label{subsec:pleadmf}
 
\noindent {\bf\boldmath $p$-leader structure functions.} \quad The $p$-leader multifractal formalism relies on the  \emph{$p$-leader structure function} $S^{(p)}_ \ell(j,q)$, which is defined as the sample moments of the $p$-leaders
\begin{equation}
\label{equ-WSF}
\forall q \in \RR, \qquad S^{(p)}_ \ell(j,q) = 2^{dj} \displaystyle \sum_{ k }  ( \pL_{j,k})^q.   
\end{equation}

\noindent {\bf\boldmath $p$-leader scaling function.} \quad The $p$-leader scaling function is defined as
\begin{equation} 
\label{qscaldef}  
\forall q \in \R, \hspace{6mm} 
  \zeta^{(p)}(q) =   \displaystyle\liminf_{j \rightarrow -\infty} \;\; \frac{\log \left( S^{(p)}_ \ell(j,q) \right) }{\log (2^{j})} \; . 
  \end{equation}
The function $ \zeta^{(p)}(q) $ does not depend on the wavelet basis when the mother wavelets are $C^{\infty}$ (e.g., in the Schwarz class), cf. \cite{JaffMel, JaffToul, Jaffard2004}. \\

\noindent {\bf\boldmath $p$-leader multifractal formalism.} \quad A multifractal formalism can be defined via the Legendre transform of the mapping $q \rightarrow  \zeta^{(p)}(q) $
\begin{equation}
\label{equ:mfhp}
{\cal L}^{(p)} (h) :=  \inf_{q \in \real} (d+qh - \zeta^{(p)}(q)). 
\end{equation}
It provides practitioners with an upper bound for ${\cal D}^{(p)} $:  If $\eta (p) >0$, then 
\begin{equation}
\label{equ-tlL}
 \forall h , \qquad {\cal D}^{(p)}(h)   \leq {\cal L}^{(p)} (h)
\end{equation}  
as a particular occurrence of the general principles that lead to (\ref{equ-tl}), see 
 \cite{JaffMel, JaffToul, JAFFARD:2010:A}.

A key property of  (\ref{equ-tlL}) is that it yields a nontrivial upper bound also for the decreasing part of the spectrum (a property that did not hold for formulas based on wavelet coefficients instead of wavelet leaders even in the setting supplied by  the H\"older exponent \cite{Jaffard2004}). 
This upper bound turns out to be sharp for many models:  this has been  verified   either theoretically (e.g., for fBm or random wavelet series \cite{JAFFARD:2010:A}), or numerically (e.g., for multifractal random walks or L\'evy processes, see \cite{Jaffard2004}).  
The generic validity of the multifractal formalism has been proven in \cite{Fraysse2}. 

The $p$-leader multifractal formalism and its estimation performance are illustrated and studied in Section~\ref{sec-num}.

\subsection{ $p$-leader Cumulants}  
 \label{subsec:cumul}

In \cite{Castaing1993}, the use of the cumulants of the log of $p$-leaders has been motivated by the fact that they enable to capture, with a small number of coefficients, most of the information contained in the spectra ${\cal L}^{(p)} (h)$.
Let $C^{(p)}_m(j) $ denote the $m$-th order cumulant of the random variables $\log (\pL_\la)$.
Assuming that the moments of order $q$ of the $p$-leader $ \pL_\la $ exist and
 \[ \EX{(\pL_\la)^q}  =2^{j \zeta^{(p)}(q)}  \EX{(\pL_{\la_0})^q } , \]  
where $\la_0$ is the unit cube $[0, 1]^d$, one obtains  that
\begin{equation}
\label{equ:cum1}
\log \left(\EX{(\pL_\la)^q} \right) =  \log \left(\EX{(\pL_{\la_0})^q} \right) +  \zeta^{(p)}(q) \log (2^{j})
\end{equation}
and, for $q$ close to $0$
\begin{equation}
\label{equ:cum2}
\log \left(\EX{(\pL_\la)^q} \right)  =  \log \left(\EX{e^{q \log \pL_\la}}\right) =  \sum_{m\geq 1} C_m^{(p)}(j) \frac{q^m}{ m!} . 
\end{equation}
Comparison of \eqref{equ:cum1} and \eqref{equ:cum2} yields that the cumulants $ C_m^{(p)}(j) $ are necessarily of the form
\begin{equation}
\label{equ-leadercum}
 C^{(p)}_m(j) = C^{(p,0)}_m + c^{(p)}_m \log (2^{j} ) ,
\end{equation}
and that $\zeta^{(p)}(q)$ can be expanded around $0$ as
\begin{equation}
\label{equ-leadercumzq}
 \zeta^{(p)}(q)  =\sum_{m\geq 1} c^{(p)}_m \frac{q^m}{ m!}. 
\end{equation}
The concavity of $ \zeta^{(p)}$ implies that  $c^{(p)}_2 \leq 0$. 
Note that, even if the moment of order $q$ of $\pL_\la $ is not finite, the cumulant of order $m$ of $\log (\pL_\la) $ is likely to be finite and \eqref{equ-leadercum} above is also likely to hold. 

Relation \eqref{equ-leadercum} provides practitioners with a \emph{direct} way to estimate the coefficients $c^{(p)}_m$ in the polynomial expansion \eqref{equ-leadercumzq}, termed the \emph{log-cumulants}, by means of linear regressions of $ C^{(p)}_m(j)$ versus $\log  (2^{j} ) $ \cite{Wendt:2007:E}. 
Moreover,  on condition that $c_2^{(p)} \neq 0$, the polynomial expansion \eqref{equ-leadercumzq} can be translated into an expansion of ${\cal L}^{(p)}(h)$ around its maximum via the Legendre transform
\begin{equation}
{\cal L}^{(p)}(h) = d + \sum_{m\geq 2} \frac{\mathcal{C}_m}{m!}  \left(  \frac{h-{c^{(p)}_1}}{c^{(p)}_2}\right)^m
\end{equation}
with 
$\mathcal{C}_2=c_2^{(p)}$, 
$\mathcal{C}_3=-{c^{(p)}_3}$, 
$\mathcal{C}_4=-{c^{(p)}_4}+3 \frac{{c^{(p)}_3}^2}{c^{(p)}_2}$, etc., see  \cite{WENDT:2009:C}.

\subsection{Estimation and finite size effects}
\label{sec:finite}

In the computation of the multifractal formalism for data in applications (in contrast with the theoretical analysis of functions) one is confronted with the fact that only values of $X$ sampled at finite resolution are available (in contrast with values for an interval in $\RR^d$). As a consequence, the infinite sum in the theoretical definition of $p$-leaders, cf.  \eqref{pleaders}, is truncated at the finest available scale induced by the resolution of the data. Moreover, the $\lim\inf$ in the definition of the scaling function \eqref{qscaldef} cannot be evaluated.
 
Our goal is here to compute the equivalents of \eqref{equ-WSF}, \eqref{qscaldef} and \eqref{equ-leadercum} for such finite resolution $p$-leaders and to use these expressions for identifying estimators for the $p$-leader scaling function $\zeta^{(p)}(q)$, the log-cumulants $c_m^{(p)}$ and the $p$-spectrum ${\cal L}^{(p)} (h)$. 
As a constructive model, we make use of a binomial deterministic wavelet cascade. For simplicity, we consider the univariate case $d=1$. 
In Section \ref{sec-num}, numerical results are provided that demonstrate that this simple deterministic model yields estimators with excellent performance for large classes of stochastic multifractal processes.

\subsubsection{ $p$-leaders for deterministic wavelet cascades}

Let  $\bar j < 0$ denote the finest and $j=0$ the coarsest available scales, respectively, i.e., $\bar j\leq j \leq 0$. Let $\w,\ww>0$ and, by convention, let $c_{j=0,k=1}=1$. The coefficients $c_{\la}$ of the deterministic wavelet cascades  are constructed by the following iterative rule \cite{m74}. 
For $j=0,-1,\ldots,\bar j+1$, for each coefficient $c_{j,k}$ at scale $j$ two children coefficients are obtained at scale $j-1$ by multiplication with $\w$ and $\ww$, $c_{j-1,2k-1}=\w c_{j,k}$ and $c_{j-1,2k}=\ww c_{j,k}$.
At a given scale $j < 0$, there are hence $2^{-j}$ coefficients, 
taking the values
$$
c_{j,k}\in\{\w^{-n}\ww^{-j-n},\;n=0,\dots,-j\}.
$$
The wavelet structure function \eqref{equ-WSF} is therefore given by
$$
S_c(j,q)=2^{j}\sum_{k=1}^{2^{-j}}(\cjk)^q = 2^{j}(\w^q+\ww^q)^{-j}=\left(\frac{\w^q+\ww^q}{2}\right)^{-j}=:2^{j\eta(q)}
$$
from which we identify the wavelet scaling function \eqref{defscalond} of the cascade
$$
\eta(q)=1-\log_2(\w^q+\ww^q).
$$
Let $p>0$ be such that $\eta(p)> 0$. 
The restricted $p$-leaders, denoted $\cpL_{\lambda,\bar j}$, are defined by replacing $3\lambda$ with $\lambda$ in \eqref{pleaders} and the corresponding structure function, denoted $\bar S^{(p)}_{\bar\ell}(j,q;\bar j)$,  by replacing $\pL_\lambda$ with $\cpL_{\lambda,\bar j}$ in \eqref{equ-WSF}.  
It can be shown that structure functions with restricted $p$-leaders yield quantities  equivalent to \eqref{equ-WSF} so that the corresponding scaling functions \eqref{qscaldef} coincide, see \cite{Bergou}.
Then
$$
\cpL_{\lambda,\bar j}:=\left(\sum_{\lambda'\subset\lambda}|c_{\lambda'}|^p2^{j'-j}\right)^{\frac{1}{p}}
    =c_{\lambda}\left(\sum_{l=0}^{j-\bar j}\left(\frac{\w^p+\ww^p}{2}\right)^{l}\right)^{\frac{1}{p}}
\!\!=c_{\lambda}\left(\sum_{l=0}^{j-\bar j}2^{-\eta(p)l}\right)^{\frac{1}{p}}.
$$
For infinite resolution $\bar j\to-\infty$, this expression boils down  to $\cpL_{\lambda,-\infty}=c_\la\left(\frac{1}{1-2^{-\eta(p)}}\right)^\frac{1}{p}$. In this case, the structure function reads
$$
\bar S^{(p)}_{\bar\ell}(j,q;-\infty)=2^{j}\sum_{k=1}^{2^{-j}}(\cjk)^q\left(\frac{1}{1-2^{-\eta(p)}}\right)^\frac{q}{p} = 2^{j\eta(q)}\left(\frac{1}{1-2^{-\eta(p)}}\right)^\frac{q}{p}
$$
from which we identify the $p$-leader scaling function  \eqref{qscaldef}  of the cascade
$$
\zeta^{(p)}(q)\equiv\eta(q).
$$ 
However, for finite resolution $\bar j>-\infty$ we have $\cpL_{\lambda,\bar j}=c_\la\left(\frac{1-2^{-( j-\bar j+1)\eta(p)}}{1-2^{-\eta(p)}}\right)^\frac{1}{p}$ and hence
\begin{equation}
\label{equ:sjqfinite}
\bar S^{(p)}_ {\bar\ell}(j,q;\bar j)= 2^{j\zeta^{(p)}(q)}\left(\frac{1-2^{-( j-\bar j+1)\eta(p)}}{1-2^{-\eta(p)}}\right)^\frac{q}{p}.
\end{equation}
Similarly, if $\bar C^{(p)}_m(j)$ denotes the $m$-th order sample cumulant of $\log (\cpL_{\la,\bar j})=\log(c_\la)+\frac{1}{p}\log\left(\frac{1-2^{-( j-\bar j+1)\eta(p)}}{1-2^{-\eta(p)}}\right)$, then
\begin{equation}
\label{equ:c1finite}
\bar C^{(p)}_1(j)=  C^{(p,0)}_1 + c_1^{(p)} \log (2^{j} ) + \frac{1}{p} \log\left(\frac{1-2^{-( j- \bar j+1)\eta(p)}}{1-2^{-\eta(p)}}\right) 
\end{equation}
while
$\bar C^{(p)}_m(j)\equiv C^{(p)}_m(j)= C^{(p,0)}_m + c_m^{(p)} \log (2^{j} )$ for $m\geq2$. Comparing  \eqref{equ:sjqfinite} and \eqref{equ:c1finite} with \eqref{qscaldef} and \eqref{equ-leadercum} shows that the truncation of the definition of $p$-leaders \eqref{pleaders} at finite scale induces an additional scale-dependent term, parametrized by $\bar j$ and $\eta(p)$, in the $p$-leaders structure functions and first log-cumulant. 

\subsubsection{Estimation}
Assuming that $\eta(p)$ is known, the expressions \eqref{equ:sjqfinite} and \eqref{equ:c1finite} enable us to define estimators for $\zeta^{(p)}(q)$ and $c^{(p)}_m$ as linear regressions in logarithmic coordinates
\begin{align}
\label{equ:estzq}
\hat\zeta^{(p)}(q) &=   \sum_{j=j_1}^{j_2} w_j\left( \log_2 \left( S^{(p)}_ \ell(j,q) \right)-\frac{q}{p}\log_2\left(1-2^{-( j-\bar j+1)\eta(p)}\right)\right) \\
\label{equ:estcm1}
\hat c^{(p)}_1 &=\frac{1}{\log (2)}\sum_{j=j_1}^{j_2} w_j\left(C^{(p)}_1(j) - \frac{1}{p} \log\left(\frac{1-2^{-( j- \bar j+1)\eta(p)}}{1-2^{-\eta(p)}}\right)\right)\\
\label{equ:estcm2}
\hat c^{(p)}_m &=\frac{1}{\log (2)}\sum_{j=j_1}^{j_2} w_jC^{(p)}_m(j),\quad\quad m\geq2.
\end{align}
Estimators for the Legendre spectra can be defined in a similar way for the parametric development $\mathcal{L}(q),\;h(q)$ described in \cite{chhabra89}.
Note that the expressions \eqref{equ:estzq} and \eqref{equ:estcm1} differ from those employed in the standard wavelet leader setting \cite{Jaffard2004,JAFFARD:2006:A,Wendt:2007:E,WENDT:2009:C,Abel} due to the additional scale-dependent terms identified in \eqref{equ:sjqfinite} and \eqref{equ:c1finite}.
However, since $\eta(p)\rightarrow\infty$ when $p\rightarrow\infty$, these expressions are equivalent to the wavelet leader ones.
In practice, the unknown function $\eta(p)$ in \eqref{equ:estzq} and \eqref{equ:estcm1} is replaced with the estimate
\begin{equation}
\label{equ:LFeta}
\hat\eta(p) =   \sum_{j=j_1}^{j_2} w_j \log_2 \left( S_c(j,p) \right).
\end{equation}
In the above expressions, $j_1$ and $j_2$ are the finest and coarsest scales, respectively, over which the estimation is performed.
The linear regression weights $w_j$ have to satisfy the constraints
$ \sum_{j_1}^{j_2} j w_{j}$   $\equiv  1$ and  $\sum_{j_1}^{j_2} w_{j}  \equiv  0$
and can be selected to reflect the confidence granted to each $\log_2 (S_\pL(j,q))$ (or $C_m^{(p)}(j)$), see e.g.  \cite{aftv00}. 
In the numerical experiments reported in Section \ref{sec-num}, following \cite{aftv00}, we perform weighted linear fits; alternative choices have been reported in \cite{Wendt:2007:E}.

\section{Multifractal Detrended Fluctuation Analysis}
\label{sec-dfa}

\subsection{Detrended Fluctuation Analysis for the estimation of the self-similarity parameter}

Detrended Fluctuation Analysis (DFA) was originally proposed for the estimation of the self-similarity parameter $H$ for fractional Brownian motion (fBm) \cite{peng1994mosaic,kantelhardt2001detecting}. 
Therefore, DFA has essentially been stated and studied in the univariate setting, $d=1$. 
In essence, it relies on the observation that the H\"older exponent of a sample path of fBm is  a constant function $h(t) \equiv H$. 
From the definition of the H\"older exponent (cf., \eqref{equ:tpe}) a \emph{natural} multiresolution quantity $T(a,t)$ emerged
\begin{equation}
\label{equ:tx:dfa}
 \TDFA(a,t) = |X(t) - P_{t,a,N_P}(t)|
\end{equation}
where $P_{t,a,N_P}$ is a polynomial of degree $N_P$ that is obtained by \emph{local fit to the data} \review{in a window of size $a$} (note that, alternatively, the use of a moving average model for $P_{t,a,N_P}$ has been proposed more recently, cf., e.g., \cite{xu2005quantifying}). 
From discrete versions $\TDFA(a,k)$ of these multiresolution quantities, structure functions 
$$
S_d(a,2)  = \frac{1}{n_a } \sum_k \TDFA(a,k)^2
$$  
are computed, and the 
parameter $H$ is classically estimated by linear regression of $\log_2 S_d(a,2) $ versus $\log_2 a$.
The performance of this procedure for the estimation of $H$ has been studied and compared to that of others, notably those based on wavelets,  in various contributions (cf. e.g., \cite{torresabry2004,oswiecimka2006wavelet}). 

\subsection{Multifractal extension}

\subsubsection{Natural extension}
To extend DFA to multifractal analysis, a straightforward choice could have been to construct a multifractal formalism based on $ \TDFA(a,t)$ \review{(cf. \cite{weber2001spectra} for a preliminary attempt)}:
$$
S_d(a,q) = \frac{1}{n_a} \sum_k \TDFA(a,k)^q \simeq C_q a^{\zeta_d(q)}, \, a \rightarrow 0.
$$
However, it is now well understood that the Legendre transform of $\zeta_d(q)$ yields a poor upper bound for the multifractal spectrum ${\cal D}(h)$.
One reason for this is that the values of $ \TDFA $  concentrate around $0$ and can thus not be raised to negative powers (note that this is also the case for increments or wavelet coefficients, leading to similarly poor estimates for ${\cal D}(h)$).
From a theoretical point of view, this poor performance can be understood in the light of the fact that the definition of the H\"older exponent \eqref{equ:tpe} leads to the incorrect intuition that it must be tied to increments, $X(t+a)-X(t)$, or ``improved'' increments, $\TDFA $,  as multiresolution quantities.
However, recent theoretical contributions (cf., e.g., \cite{Jaffard2004,JAFFARD:2006:A,Wendt:2007:E,barralgoncalves,Abel}) show that correct multiresolution quantities associated with the H\"older exponent are the \emph{oscillations}
$$ {\cal O}(a,x) = \sup_{u, v \in B(x,a) }| X(u) - X(v)|.$$
A relevant multifractal formalism can be constructed for ${\cal O}(a,x)$ as long as the smooth parts of the data have H\"older exponents strictly smaller than $1$. 
Otherwise, higher order oscillations must be used.
Note that in a wavelet framework, wavelet leaders play for wavelet coefficients the same role as oscillations play for increments.

\subsubsection{ $L^2$-norm formulation}

The use of oscillations $ {\cal O}(a,x)$ instead of $\TDFA$ would thus have been a natural way to extend DFA to multifractal analysis.
However, the multifractal extension proposed by Kantelhardt et al. followed  a different path: In the seminal contribution \cite{kantelhardt2002multifractal}, they  proposed the following original multiresolution quantity
\begin{equation}
\label{equ:mradfa}
\TMFD(a,k) =  \left( {1\over a} \sum_{i=1}^{a} |X(ak+i) - P_{k,a, N_P}(i) |^2 \right)^\frac{1}{2}, k =1, \ldots, \frac{n}{a}
 \end{equation}
and constructed a multifractal formalism with $\TMFD(a,k)$, the so-called \emph{multifractal detrended fluctuation analysis} (MFDFA).
Here $n $ denotes the number of available samples and $P_{k,a, N_P}$ is the same polynomial as in \eqref{equ:tx:dfa}. 
Defining structure functions
\review{$$
S_{mfd}(a,q)=  \frac{a}{n}  \sum_k^{n/a} \TMFD(a,k)^q =  \frac{a}{n}  \sum_k^{n/a} \left( {1\over a} \sum_{i=1}^{a} |X(ak+i) - P_{k,a, N_P}(i) |^2 \right)^{\frac{q}{2}}
$$}
and the corresponding scaling function
$$ 
\zeta_{mfd}(q) =   \displaystyle\liminf_{j \rightarrow -\infty} \;\; \frac{\log \left( S_{mfd}(a,q)   \right) }{\log (a)}
 $$
yields, via a Legendre transform, the MFDFA multifractal formalism
 $$
 \LMFD(h) = \inf_{q \in \real} (d+qh - \zeta_{mfd}(q)) \geq {\cal D}(h).  
 $$
It was numerically compared against the wavelet leader and WTMM multifractal formalisms for multiplicative cascade-type multifractal processes, see \cite{kantelhardt2002multifractal} for the original contribution.
Despite the fact that there has been no theoretical motivation for the use of \eqref{equ:mradfa}, the MFDFA formalism was found to perform very satisfactorily and is commonly used in applications (cf., e.g., \cite{galaska2008comparison,lim2007multifractal,Wang2011goldmarket,shang2008detecting,Yuan2009,Wang2012}), thus empirically justifying a posteriori the choice of $ \TMFD$ as a relevant multiresolution quantity for multifractal analysis.

\subsection{ $p$-exponents, $p$-leaders and MFDFA}
\label{sec:PLmfdfa}
Comparing the definition of the multiresolution quantity $T^{(p)}$ in \eqref{equ:tpa} for the $p$-exponent to that of $\TMFD$ in \eqref{equ:mradfa} for MFDFA clearly shows that the latter mimics the former, for $p=2$, with a discretized setting of the continuous integral in \eqref{equ:tpa}, and with the \emph{theoretical} Taylor polynomial $P_{x_0}$ replaced by a \emph{data-driven locally adjusted} polynomial $P_{k,a, N_P}$.
Therefore, MFDFA can a posteriori be interpreted, and theoretically grounded, in the framework of $p$-exponent analysis:
While MFDFA was originally associated to the analysis of local regularity as measured by the H\"older exponent, the present contribution thus shows that  $\TMFD$ \emph{must be related to a $p$-exponent based characterization of local regularity with $p=2$ and not to the H\"older exponent.}

The following paragraphs provide a detailed theoretical and practical comparison of the MFDFA and $p$-leader frameworks. 
\vskip2mm

\noindent {\bf\boldmath Choosing $p$.} \quad In the MFDFA method, $p$ is arbitrarily set to $2$, while $p$-exponents in  \eqref{equ:tpa} are theoretically defined for $p \in (0, p_0)$. 
There are practical benefits stemming from varying $p$, as summarized in Section~\ref{sec:pproperties} and discussed in detail in the companion paper \cite[Section 4]{PART1}:
Notably, the analysis requires to choose $p$ such that $ \eta (p) >0$ ; moreover,
the use of various values $p < p_0$ may help practitioners to reveal the fine local singularity structure of data (see Section \ref{sec:lacunary} for a numerical illustration). 
\vskip2mm

\noindent {\bf Local regularity and integration.} \quad In numerous real-world data, notably in biomedical applications (Heart Rate variability, fMRI), it is observed that $\Hmin < 0$ (cf., e.g., \cite{MandMemor}). 
This explains why MFDFA procedures, as detailed, e.g., in \cite{kantelhardt2002multifractal,hardstone2012detrended,ihlen2012introduction}, always perform an integration of the data, $ X(t) \rightarrow \int^t  X(s) ds$, as a preliminary step (hence implicitly assuming that for the integrated data, $\Hmin>0$).
However, as long as $\eta(p) > 0 $ and $h > -1/p $, this preprocessing of data is not needed and may even constitute a drawback if data contain oscillating singularities, cf.  \cite[Section 4]{PART1}. 

\vskip2mm

\noindent {\bf Time domain versus wavelet domain.} \quad The MFDFA method relies on a time domain implementation, i.e., $\TMFD$ is computed directly from the sampled time series $X$, while the $p$-leader formalism relies on wavelet projections. This has fundamental theoretical and practical implications which are detailed in the next paragraphs.
\vskip2mm

\noindent {\bf Local regularity, Taylor polynomial and vanishing moments.} \quad 
A key property of the wavelet characterization is that it does not require the knowledge of the Taylor polynomial (cf. \cite{Abel}) whereas  MFDFA require some (heuristic) estimation of this polynomial \review{(cf. a contrario \cite{gu2010detrending,schumann2011multifractal}, where a variation of MFDFA that attempts to avoid the estimation of the polynomial is devised)}.
Despite its being inspired from the definition of the H\"older exponent, the polynomial $P_{t,a,N_P}(t)$ is in practice obtained as the \emph{best fit} of $X(u)$ for $u \in [t-a/2, t+a/2]$ for a priori chosen and fixed degree $N_P$. The polynomial $P_{t,a,N_P}$ must be computed for  each time position $t$ and analysis scale $a$. It thus depends on $a$, while it should \emph{not} depend on $a$ in theory, cf., \eqref{equ:tpa}. Moreover, the order of the polynomial in  \eqref{equ:tpa}  can depend on the time position $t$ (cf., \cite{Abel}) while the order $N_P$ in the MFDFA method is fixed.

Another key difference in the use of $p$-leaders $\pL$ and $\TMFD$ lies in the fact that the former requires the choice of the number of vanishing moments $N_\psi $ of the mother wavelet $\psi_0$, while the latter implies the choice of the degree $N_P$ of the polynomial $P_{t,a,N_P}$. 
Often, the parameters $N_\psi $ and $N_P$ are regarded as equivalent, yet this is only partially correct since the choices of $N_\psi $ and $N_P$ are framed by different theoretical constraints:
In order to recover the local power-law behavior $ T^{(p)}_X(a,t) \sim C a^{h}, a \rightarrow 0$ for an isolated singularity with regularity exponent $h$, $ N_\psi $ must be larger than $h$ (and can be set globally to be larger than  the largest regularity exponent of $X$, cf., \cite[Section 3.1]{PART1}, \cite{Jaffard2004}), while it is required that $N_P<h$ by definition of the H\"older exponent. 
\vskip2mm

\noindent {\bf Robustness to trends and finite size effects.} \quad 
Practically, the choice of $N_\psi $ and $N_P$ is often thought of in terms of robustness to trends (cf. e.g., \cite{va99,peng1994mosaic,torresabry2004,HuDFA2001,Nagarajan2005,ChenDFA2001,Bashan2008}).

Assuming that the data to analyze are actually corrupted by deterministic trends $Z$ seen as noise, $X + Z$, it has been documented that increasing the number of vanishing moments of the mother wavelet $N_\psi $ diminishes the impact of the additive trends to the estimation of the scaling exponent \cite{va99,torresabry2004}: The possibility of varying $N_\psi$ brings theoretical and practical robustness to analysis and estimation.
The practical price paid for increasing $N_\psi $ consists of a larger number of wavelet coefficients being polluted by border effects (finite size effect) and can also decrease the coarsest scale at which data can be analyzed. 
This price turns out to be very low since coarse scales contain only few coefficients and hence have very little impact on estimated scaling exponents. Note that the finest scale that is available using $p$-leaders is independent of $N_\psi$ and directly and only determined by the resolution of the data (i.e., if the data is sampled at $\Delta t$, then the the finest scale contains coefficients at rate $2\Delta t$).

In the MFDFA method, increasing $N_P$ also brings some form of robustness to additive noise, however, this far more depends on the nature of the trend \cite{HuDFA2001,Nagarajan2005,ChenDFA2001,Bashan2008} (see also the numerical illustrations in Section~\ref{sec:trend}).
Yet, increasing $N_P$ has a more drastic practical consequence: The finest scales of data cannot be used practically since the \emph{best fit} of $X(u)$ for $u \in [t-a/2, t+a/2]$ cannot be achieved until the number of samples actually available in $[t-a/2, t+a/2]$ is substantially larger than $N_P$: For a polynomial of order $N_P$, one needs at least $N_P+2$ data points and consequently, $\Delta j=\lceil\log_2(N_P+2)\rceil-1$ fine scales that can be analyzed with $p$-leaders are not accessible for MFDFA. 
Losing fine scales is problematic from a statistical estimation point of view (since fine scales contain many coefficients) as well as from a multifractal analysis point of view (since it theoretically requires to estimate the scaling exponents in the limit of fine scales $a \rightarrow 0$). 

In summary, the choice of $N_\psi $ is theoretically better grounded and practically much easier and less critical than that of $N_P$.
The polynomial subtraction entering MFDFA is thus a more intricate issue than it may seem, with little theoretical guidelines. 
\vskip2mm

\noindent {\bf Extension to higher dimension.} Theoretically, MFDFA
can be extended straightforwardly to higher dimensions, $d >1$.
In practice, the computation of local best fit polynomials becomes a real issue in higher dimension, and there are few attempts to extend MFDFA to dimension $d=2$ \cite{Gu2006,soares2009towards,wang2014local,shi2015new}. 
Along the same line, the use of polynomials of degree 2 or higher is uneasy (higher order polynomial fitting and multivariate integral numerical evaluation). In contrast, the $p$-leader formalism extends without difficulties to higher dimensions given that higher dimensional discrete wavelet transforms are readily obtained by tensor product of 1D-mother wavelets (cf., e.g., \cite{Antoine2004}). 
Also note that the MFDFA method has larger computational complexity (of order $\mathcal{O}(n\log n)$, where $n$ is the total number of samples) than the $p$-leader formalism (of order $\mathcal{O}(n)$), which quickly becomes an issue in higher dimension.
\vskip2mm

\noindent {\bf Conclusions.} \quad This connection of MFDFA  to $p$-leaders provides a theoretical grounding for the choice of $\TMFD$, which otherwise appears as a relevant, yet ad-hoc, intuition used to construct a multifractal formalism.
It also makes explicit that MFDFA measures local regularity via the $2$-exponent and \emph{not} the H\"older exponent. 
The $p$-leader formalism can thus be read as an extension (different $p$) and wavelet-based reformulation of MFDFA. 
The $p$-leader and MFDFA formalisms are compared in terms of estimation performance and robustness to trends in Section~\ref{sec-estim} below. 

It is also of interest to note that $p$-exponents and $p$-leaders, on one side, and MFDFA, on other side, were proposed independently and in different fields: 
In the Mathematics literature for the former around year 2005 \cite{JaffMel,JaffToul}, in the Physics literature for the latter around year 2002 \cite{kantelhardt2002multifractal}.
This is, to the best of our knowledge, the first time that these two notions are connected, related and compared.

\section{Illustrations and estimation performance}
\label{sec-num}

We numerically illustrate and validate the proposed $p$-leaders multifractal formalism and compare it against the leader and MFDFA formalisms. 
To this end, the formalisms are applied to independent realizations of synthetic random processes with prescribed multifractal $p$-spectra, and their respective estimation performances are studied in detail and compared.
For the sake of exhausitivity in comparisons, processes whose $p$-spectra do not depend on $p$ as well as processes whose spectra vary with $p$ are used. 
Also, situations in which $p$-exponents all collapse with the H\"older exponents are considered as well as examples where this is not the case. 
The $p$-leader multifractal formalism is furthermore numerically validated in higher dimensions by application to synthetic (2D) images with prescribed multifractal $p$-spectra. 

The WTMM formalism has already been compared independently against the wavelet-leader formalism \cite{JAFFARD:2010:A,JAFFARD:2006:A} and MFDFA \cite{kantelhardt2002multifractal,oswikecimka2012effect} and is thus not re-included in this study. 

Results illustrate the benefit of the extra flexibility of varying $p$ in the $p$-leader multifractal formalism, both in terms of estimation performance and for evidencing data whose multifractal $p$-spectra are not $p$-invariant.

Sample paths of all processes were numerically simulated by {\sc Matlab} routines implemented by ourselves, available upon request. 

\subsection{ Random processes with prescribed multifractal $p$-spectra}
\label{sec-proc}

\subsubsection{Fractionally differentiated Multifractal Random Walk}

The multifractal random walk (MRW) has been introduced in \cite{bdm01} and is defined as
$$
X(k)=\sum_{k=1}^{n} G_H(k) \exp({\omega(k)})
$$
where $G_H(k)$ are increments of fractional Brownian motion with parameter $H>1/2$ and $\omega$ is a Gaussian process that is independent of $G_H$ and has covariance  
$
\textnormal{cov}(\omega(k_1),\omega(k_2))=\lambda \ln\left(\frac{L}{|k_1-k_2|+1}\right)
$
when $|k_1-k_2|<L$, and $0$ otherwise, with $\lambda> 0$. 
By construction, MRW has stationary increments and mimics the multifractal properties of Mandelbrot's multiplicative cascades \cite{m74}. 
It has been chosen here as an easy-to-simulate member of this widely used class of multifractal processes. \\
By means of fractional integration of negative order $s<0$ of $X$ (in practice fractional differentiation of positive order $\nu=-s$, cf., \cite[Definition 3]{PART1}), we obtain sample paths $X^{(\nu)}$ with different values for the critical Lebesgue index $p_0$.
MRW contains only canonical singularities (cf. \cite[Definition 4]{PART1}), hence fractional differentiation results in a pure shift of their multifractal $p$-spectra by $\nu$ to smaller values of $h$. 
Moreover, since canonical singularities are $p$-invariant (cf., \cite[Theorem 2]{PART1}), the multifractal $p$-spectra of $X^{(\nu)}$ collapse for all $ p\leq p_0$ and are given by
\begin{equation}
\label{equ:mrwfdDH}
\mathcal{D}^{(p)}_\nu(h)\equiv \mathcal{D}_\nu(h)=
\begin{cases}
1+\frac{c_2}{2}\left(\frac{h-c_{1,\nu}}{c_2}\right)^2& \mbox{for } h_p\in\Big[c_{1,\nu}-\sqrt{-2c_2},\;c_{1,\nu}+\sqrt{-2c_2}\Big]\\
-\infty & \mbox{otherwise}
\end{cases}
\end{equation}
where $c_{1,\nu}=H+\lambda^2/2-\nu$, $c_2=-\lambda^2$. 
Furthermore, $c_m\equiv 0$ for all $m\geq3$.
The wavelet scaling function is given by
$$
\eta(p)=
\begin{cases}
(H+\lambda^2/2-\nu)p+\frac{c_{2}}{2}p^2 & \mbox{for } 0\leq p\leq\sqrt{-2/c_2}\\
(H+\lambda^2/2-\nu)\sqrt{-2/c_2}-1-\sqrt{-2c_2}+c_2p & \mbox{for } p>\sqrt{-2/c_2}.
\end{cases}
$$
By elementary calculations, evaluation of condition \eqref{equ:p0} yields
\begin{equation}
\label{equ:p0mrw}
p_0=
\begin{cases}
{\infty}&\mbox{for } \nu\in\big[0,\;  H+\lambda(\frac{\lambda}{2}-\sqrt{2})\big]\\
1\big/\big(\nu-H-\lambda(\frac{\lambda}{2}-\sqrt{2})\big)& \mbox{for } \nu\in\big(H+\lambda(\frac{\lambda}{2}-\sqrt{2}),\; H+\lambda(\frac{\lambda}{2}-\frac{1}{\sqrt{2}})\big]\\
2(H+\frac{\lambda^2}{2}-\nu)/\lambda^2& \mbox{for } \nu\in\big(H+\lambda(\frac{\lambda}{2}-\frac{1}{\sqrt{2}}),\; H+\frac{\lambda^2}{2}\big]
\end{cases}
\end{equation}
and $\eta(p)<0$ for any $p>0$ if $\nu\geq H+\lambda^2/2$.

\subsubsection{Lacunary wavelet series}
\label{sec:rws}

Lacunary wavelet series (LWS) $X_{\al, \eta}$ depend on a lacunarity parameter $\eta \in (0,1)$ and a regularity  parameter  $\al \in \RR$. 
At each scale $j\leq 0$, the process  has  a fraction of exactly $2^{-\eta j}$ nonzero  wavelet coefficients on each interval $[l, l+1) $ ($l \in \ZZ$), at uniformly distributed random locations, and whose common amplitude  is $2^{\al j}$.
This process was introduced in \cite{Jaf6} when $\al >0$ (i.e., suited to H\"older exponent based multifractal analysis) and in \cite{Bandt2015} for the general case $\al \in \RR$ (thus requiring the use of $p$-exponents).
The use of LWS here is motivated by the fact that they contain singularities that are not $p$-invariant: Indeed, it is shown in \cite{Bandt2015}  that almost every point is a \emph{lacunary singularity} (cf., \cite[Definition 7]{PART1}). For $\alpha>0$, $\eta \in (0,1)$, LWS are bounded and their $p$-spectra $\forall p >0$ are given by \cite{Bandt2015}
\begin{equation}
\label{eq:DhpLWS}
\mathcal{D}^{(p)}(h)= 
\begin{cases}
\eta \displaystyle\frac{H + 1/p}{\al + 1/p}  & h\in \left[\al, \displaystyle\frac{\al}{\eta} + \left(\displaystyle\frac{1}{\eta}  -1 \right) \frac{1}{p}\right] \\ \\ 
-\infty & \mbox{otherwise}.
\end{cases}
\end{equation}

\subsubsection{Simulation setup}
\label{sec:simsetup}

We generate $N_{MC}=500$ independent realizations of $N=2^{16}$ samples each of MRW (with parameters  $H=0.72$, $\lambda=\sqrt{0.08}$ and fractional differentiation 
$\nu\in\{0,\;0.4,\;0.6,\;0.7,\;0.73\}$, yielding $p_0=\{+\infty,\;25,\;4,\;1.5,\;0.75\}$, cf., \eqref{equ:p0mrw}) and LWS (with parameters $\alpha\in\{0.2,0.3\}$, $\eta\in\{0.7,0.8\}$, $p_0=+\infty$). 
For each realization, we compute the Legendre spectra $\mathcal{L}^{(p)}(h_p)$ in \eqref{equ:mfhp} and $\LMFD(h)$ as well as the log-cumulants $c_m^{(p)}$ for \review{$m=1,2,3$} (cf., (\ref{equ-leadercum}-\ref{equ-leadercumzq})).
We adhere to the convention that the finest available dyadic scale is labelled by $\bar j=1$. 
Estimates are computed using (\ref{equ:estzq}-\ref{equ:estcm2}) for dyadic scales from $j_1=4$ to the coarsest available scale  $j_2$ ($j_2=13$ for $p$-leaders due to border effects and $j_2=15$ for MFDFA) with weighted linear regressions ($b_j\!=\!n_j $, see \cite{aftv00}). For $p$-leaders, a Daubechies' wavelet with $N_\psi=2$ vanishing moments is used, and consistently the degree of the polynomial for MFDFA is set to $N_P=1$. Furthermore, the $p$-leader estimates are calculated for the set of values $p\in\{\frac{1}{4},\;\frac{1}{2},\;1,\;2,\;4,\;5,\;8,\;10,\;+\infty\}$ (where $p=+\infty$ corresponds to wavelet leaders).

\subsection{ $p$-invariant $p$-spectra and negative regularity}

We use fractionally differentiated MRW here because it enables us to study negative regularity and to compare the respective performance for different values of $p$ since its multifractal $p$-spectra are $p$-invariant.

\subsubsection{ Numerical illustrations of multifractal $p$-spectra}
\label{sec-illust}

\begin{figure}[tb]
\centering
\includegraphics[width=\linewidth]{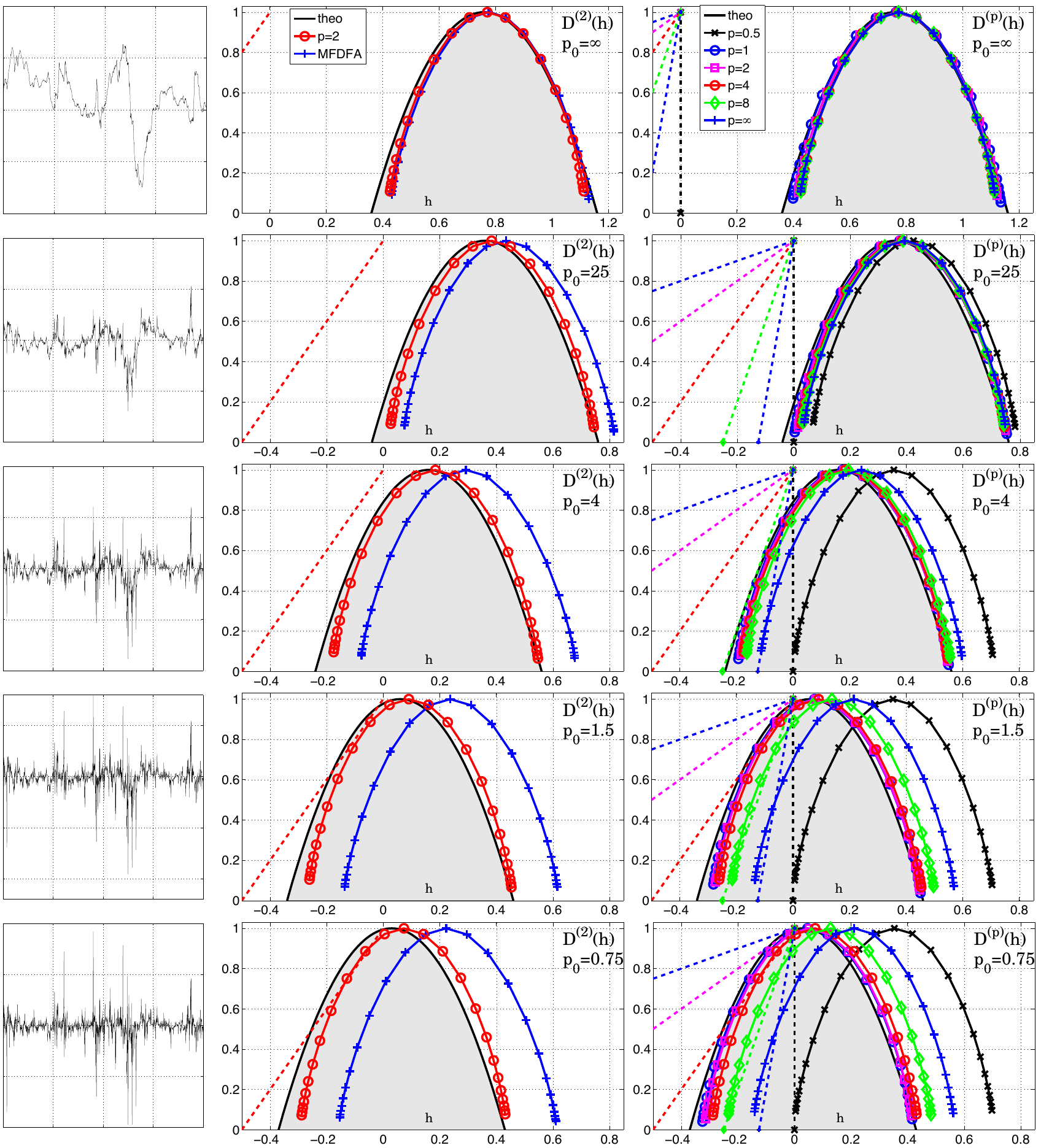}%
\caption{\label{fig:illMRW}{\bf Fractionally differentiated MRW} with $p_0=\{+\infty,\;25,\;4,\;1.5,\;0.75\}$ (from top to bottom row); the top row plots MRW without fractional differentiation: 
Single realizations (left column), theoretical spectra $\mathcal{D}(h)$ (black solid line and shaded area), estimates $\LMFD(h)$ and $\mathcal{L}^{(2)}(h)$ (center column) and estimates $\mathcal{L}^{(p)}(h)$ (right column).
The dashed lines indicate the theoretical bound \eqref{bound:dh}. 
}
\end{figure}

We illustrate and qualitatively compare the $p$-leader, leader and MFDFA multifractal formalisms, based on their Legendre spectra $\mathcal{L}^{(p)}(h)$ and $\LMFD(h)$. 
Averages over independent realizations are plotted in Fig. \ref{fig:illMRW} (colored solid lines with symbols), together with the theoretical spectra \eqref{equ:mrwfdDH} (black solid lines and shaded area) and the respective theoretical bounds  \eqref{bound:dh} for $\mathcal{D}^{(p)}(h)$ (colored dashed lines).

\vskip2mm\noindent{\bf Sample paths.}
Fig.  \ref{fig:illMRW} (left column) plots representative examples of sample paths $X^{(\nu)}$ with, from top to bottom, increasing value of $\nu$ (and, hence, decreasing regularity and $p_0$). Visual inspection of the sample paths indicates the practical benefit of the use of model processes with (potentially negative) regularity: While (positive only) H\"older regularity implies relatively smooth sample paths, the use of (negative) $p$-exponents provides practitioners with a rich set of models for applications with highly irregular sample paths, with a continuously rougher appearance as $p$-exponents take on smaller and smaller values (cf. first column of Fig.  \ref{fig:illMRW}).

\vskip2mm\noindent{\bf\boldmath $2$-leaders and MFDFA.} In Fig. \ref{fig:illMRW} (center column) Legendre spectra obtained with MFDFA and $p$-leaders with $p=2$  are plotted together with the theoretical $2$-spectra $\mathcal{D}^{(2)}(h)$. 
\begin{itemize}
\item[-] The estimates $\mathcal{L}^{(2)}(h)$ and $\LMFD(h)$ are observed to be qualitatively equivalent when $p_0$ is large. However, for negative regularity ($2<p_0\ll +\infty$), the Legendre spectra $\LMFD(h)$ obtained with MFDFA only partially capture the theoretical spectra for negative values of $h$ and appear shifted to larger values of $h$ with respect to the theoretical spectrum $\mathcal{D}^{(2)}(h)$. 
\end{itemize}

\clearpage

\begin{itemize}
\item[-]  In contrast to MFDFA, the $2$-leaders formalism clearly provides excellent estimates for $\mathcal{D}^{(2)}(h)$ for any of the processes $X^{(\nu)}$ for which $p_0\geq2$. 
\end{itemize}
We conjecture that the bias (shift towards positive values of $h$) of $\LMFD(h)$ for processes with negative regularity  is caused by finite size effects similar to those analyzed for $p$-leaders in Section \ref{sec:finite}.

\vskip2mm\noindent{\bf\boldmath $p$-leaders for $p\neq2$.} In the right column of Fig. \ref{fig:illMRW}, averages of $p$-leader estimates $\mathcal{L}^{(p)}(h)$ are plotted for several values of $p$ and compared to the theoretical $p$-spectra $\mathcal{D}_{\nu}^{(p)}(h)$ and the theoretical bounds \eqref{bound:dh}.
\begin{itemize}
\item[-] Clearly, the $p$-leader multifractal formalism provides excellent estimates for the theoretical spectra $\mathcal{D}_{\nu}^{(p)}(h)$ for $p\leq p_0$, hence validating the proposed formalism. 
\item[-] The estimates are of excellent quality also for $p<1$; notably, the choice $p=1/2<1$ enables to correctly recover the $p$-spectrum $\mathcal{D}_{\nu}^{(p)}(h)$ for $X^{(\nu)}$ with 
$p_0=0.75$ 
(Fig. \ref{fig:illMRW}, bottom row), which is not possible for $p\geq1$ (and, hence, neither with the MFDFA method). 
\item[-] When $p> p_0$, the estimates $\mathcal{L}^{(p)}(h)$ are tangent to the theoretical bounds \eqref{bound:dh}. Consequently, they are shifted to larger values of $h$ with respect to the theoretical spectrum $\mathcal{D}_{\nu}^{(p)}(h)$ and hence biased. This is visually most striking for the case $p=+\infty$ (i.e. for classical leaders associated with H\"older exponents), for which estimated spectra are constrained to positive values of $h$.
This phenomenon will be investigated in a forthcoming study (see \cite{GRETSI15-pL} for preliminary results).
\item[-] Finally, in consistency with \cite[Theorem 2]{PART1}, the spectra $\mathcal{L}^{(p)}(h)$ coincide for all $p\leq p_0$ for fractionally differentiated MRW. 
\end{itemize}

\begin{figure}[tb!]
\centering%
\includegraphics[width=\linewidth]{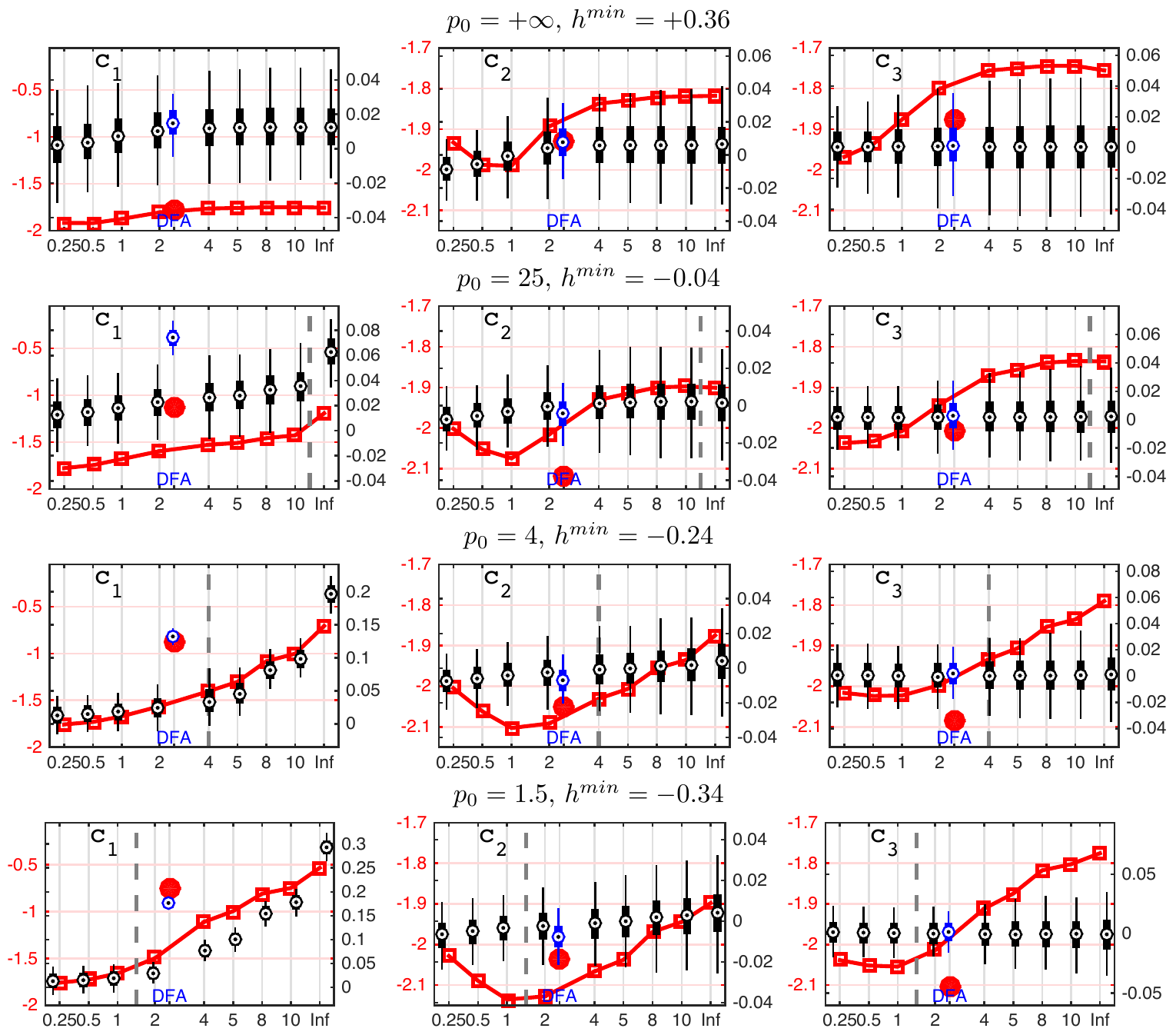}
\caption{\label{fig:estMRW}{\bf log-cumulants of fractionally differentiated
    MRW} with $p_0=\{+\infty,\;25,\;4,\;1.5\}$ (from top to bottom). \review{Red:} $\log_{10}$ of   rmse of $p$-leaders (solid lines and squares), as a function of $p$, and for   MFDFA (\review{full circles}), for $c_m$, $m=\{1,\;2,\;3\}$ (left to right   column). \review{Black and blue: boxplots of the error ($\hat c_m - c_m$) for $p$-leaders (black) and MFDFA (blue). The left and right y-scales correspond, respectively, to rmse and boxplots.  Rmse is shown in the same y-scale across the rows.} }
\end{figure}

\subsubsection{Estimation performance}
\label{sec-estim}

We proceed with a quantitative analysis of the estimation performance of the $p$-leader and MFDFA multifractal formalisms, respectively. 
To this end, we assess the estimation performance for the log-cumulants $c_m^{(p)}$ for $m=1,2,3,4$ (cf., (\ref{equ-leadercum}-\ref{equ-leadercumzq})) based on their root mean squared error (rmse), defined as
$$
\mbox{rmse}_{c_m^{(p)}}=\sqrt{\left\langle \big(\hat c_m^{(p)} -c_m^{(p)}\big)^2\right\rangle_{N_{MC}}}
$$
where $\langle\;\cdot\;\rangle_{N_{MC}}$ stands for the average over independent realizations.

Results for fractionally differentiated MRW $X^{(\nu)}$ with $\nu\in\{0,\;0.4,\;0.6,\;0.7\}$ ($p_0=\{+\infty,\;25,\;4,\;1.5\}$)
are plotted in Fig. \ref{fig:estMRW} (top to bottom rows, respectively); the logarithm ($\log_{10}$) of rmse values of $p$-leaders are plotted as a
function of $p$ (solid \review{red} lines with squares), those of MFDFA are plotted \review{with a red circle}. \review{Distributions of the estimates, after subtraction of theoretical value, are shown in black boxplots for $p$-leaders and blue boxplots for MFDFA. } The values for $p$ on the right of the
vertical red dashed lines are larger than $p_0$, $p>p_0$. Since $c_m^{(p)}\equiv c_m$, we omit the superscript $\cdot^{(p)}$ below.

\vskip2mm\noindent{\bf\boldmath Estimation of $c_1$.} Rmse values \review{and estimates} for the first log-cumulant are reported in the left column of Fig. \ref{fig:estMRW} and yield the following conclusions.
\begin{itemize}
\item[-] For $p\leq p_0$, the $p$-leader multifractal formalism systematically yields better estimation performance for small values of $p$ than for large values of $p$; the improved performance for small values of $p$ is induced by reduced standard deviations of estimates, resulting in considerable rmse gains of up to a factor $2$ for $p=1/2$ over the value of $p$ which is picked closest to  $p_0$.
\item[-] For $p>p_0$, there is a sharp increase in rmse due to a systematic bias. Indeed, $c_1$ captures the position of the maximum of the $p$-spectra, which is shifted to larger values of $h$ as compared to the theoretical spectrum when $p>p_0$, as discussed in Section \ref{sec-illust}.
\end{itemize}

\vskip2mm\noindent{\bf\boldmath Estimation of $c_m$ for $m\geq 2$.} 
The second, third and fourth column of Fig. \ref{fig:estMRW} plot rmse values \review{and estimates} for $c_2$, $c_3$ and $c_4$, respectively, and yield the following conclusions:
\begin{itemize}
\item[-] A systematic benefit in choosing small values for $p$ in the $p$-leader multifractal formalism is observed. Indeed, rmse values decrease sharply and pronouncedly for $p\leq 4$, and the smallest rmse values are obtained for small values for $p$ ($p=1/2$ or $p=1$).
\item[-] Unlike for $c_1$, choosing $p>p_0$ does not significantly alter estimation performance for $c_m$ for $m\geq 2$. 
The interested reader is referred to \cite{GRETSI15-pL} for further details. 
\end{itemize}

\vskip2mm\noindent{\bf\boldmath Comparison with MFDFA.} 
\begin{itemize}
\item[-] The MFDFA and $2$-leader formalisms have similar performance for the estimation of $c_1$ as long as $p_0=+\infty$. As soon as the support of the $p$-spectrum includes negative values for $h$ and $p_0<\infty$, MFDFA produces estimates for $c_1$ that are biased. This is consistent with the observations of Section \ref{sec-illust} where 
the $\LMFD(h)$ are found to be shifted to larger values of $h$.
\item[-] For the estimation of $c_m$ for $m\geq 2$, the MFDFA formalism yields similar performance as the $p$-leader formalism for moderately small values $p\approx 2$.
\end{itemize}
These results lead to the conclusion that, for data containing $p$-invariant singularities only, it is beneficial to choose a small value of $p$ in the analysis. Note that the $p$-leader formalism with moderately small values for $p$, e.g. $p\leq 4$, significantly outperforms the current state-of-the-art wavelet leader formalism ($p=+\infty$) which yields up to 50 percent larger rmse values. The MFDFA and $p$-leaders formalisms have comparable performance for the estimation of $c_m$ for $m\geq 2$ and also for $c_1$ as long as $p_0=+\infty$. Yet, MFDFA estimates of $c_1$ are biased when the data are characterized by negative regularity exponents.

\begin{figure}[tb]
\centering%
\includegraphics[width=\linewidth]{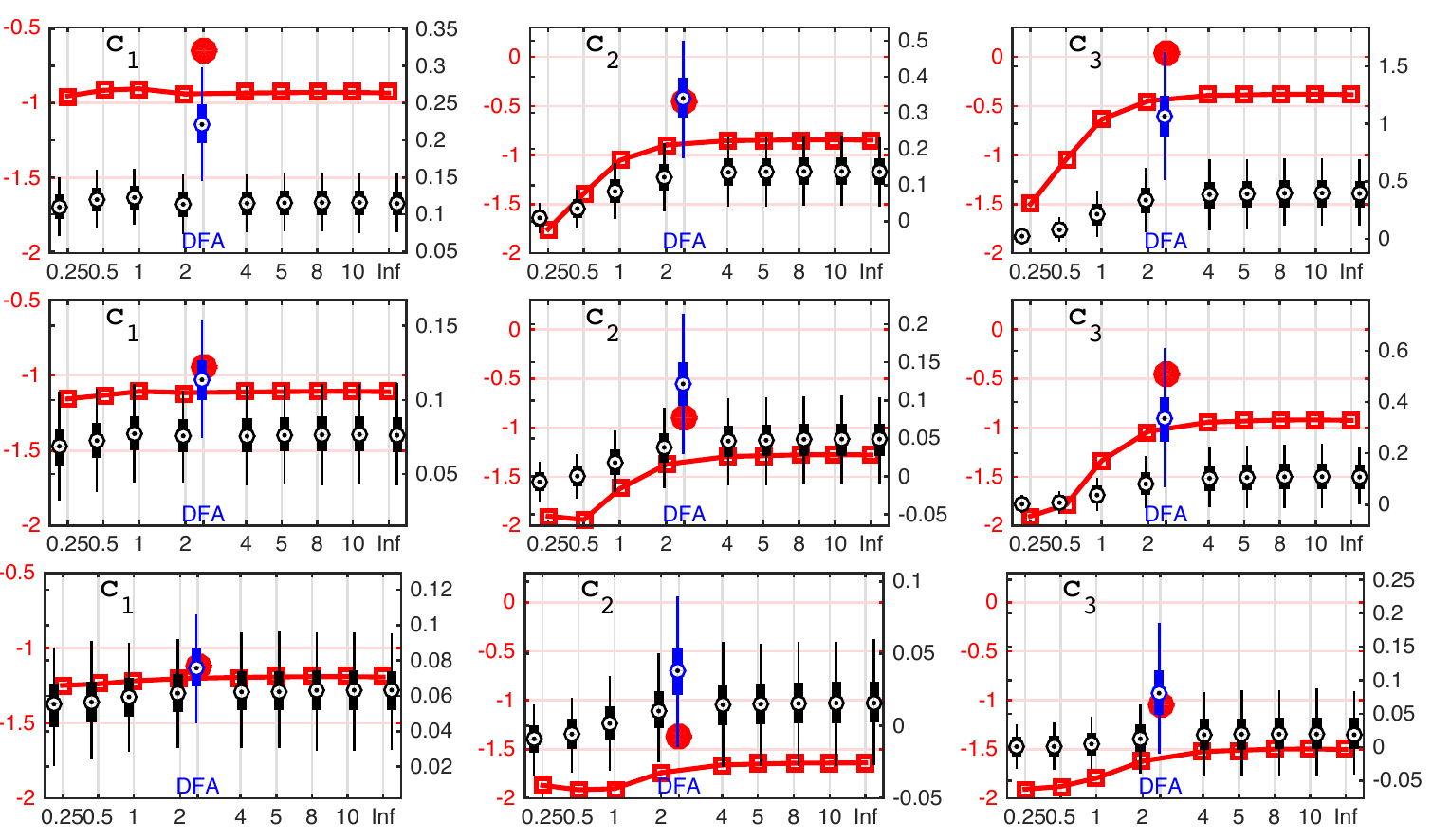}
\caption{\label{fig:estMRWcusp}{\bf MRW with non-polynomial trend} \review{($H=0.72$, $\lambda=\sqrt{0.08}$, $\nu=0$), for different detrending powers: $N_\psi=3$, $N_P=2$ (top row),  $N_\psi=4$, $N_P=3$ (center row),  $N_\psi=5$, $N_P=4$ (bottom row). Red: rmse for $p$-leaders (solid squares), as a function of $p$, and for MFDFA (full circles). Black and blue: boxplots of the error ($\hat c_m-c_m$) for $p$-leaders (black) and MFDFA (blue). The left and right y-scales correspond, respectively, to rmse and boxplots.  Rmse is shown in the same Y-scale across the rows.}}
\end{figure}

\subsubsection{Robustness against non-polynomial trend}
\label{sec:trend}
Here, we study the robustness of MFDFA and the $p$-leader formalism for estimation of the log-cumulants $c_m$ for $m=1,2,3,4$ when a $C^{\infty}$ but non-polynomial trend of the form
$$
\tau(t)=100(t+1/100)^{-1/2},\qquad t\in[0,1]
$$
is added to MRW with parameters specified as in Section \ref{sec:simsetup} and $\nu=0$.
We use Daubechies' wavelet with $N_\psi=\{3,4,5\}$ vanishing moments for $p$-leaders and, correspondingly, polynomials of degree $N_P=\{2,3,4\}$ for MFDFA. Note that with the largest choice $N_P=4$ for the degree of polynomial, the finest available scale for MFDFA is $j=4$; unlike $p$-leaders, finer scales cannot be used for estimation with MFDFA, cf., Section \ref{sec:PLmfdfa}. 
In order not to penalize MFDFA in the comparisons, we perform linear regressions from scale $j_1=4$ to the coarsest available scales ($j_2=13$ for $p$-leaders due to border effects and $j_2=15$ for MFDFA).

Rmse values \review{and estimates} for $c_m$, $m=1,\ldots,4$, are plotted in Fig. \ref{fig:estMRWcusp} for $N_\psi=\{3,4,5\}$  (top to bottom rows, respectively). 
Clearly, the estimation performance of MFDFA is severely degraded by the non-polynomial trend: Even with the polynomial of highest degree considered here ($N_P=4$), rmse values for $c_2$, $c_3$ and $c_4$ for MFDFA are up to one order of magnitude larger than in the absence of the trend (cf. Fig. \ref{fig:estMRW}, top row). 

In contrast, the rmse values for $p$-leaders with $N_\psi=4$ are very close to those obtained in the absence of the trend.
This indicates that  the wavelet transform underlying $p$-leaders is considerably more effective in removing the impact of the trend on the higher order statistics of the multiresolution quantities than the empirical polynomial-fitting procedure employed by the MFDFA method.

\begin{figure}[tb]
\centering%
\includegraphics[width=0.72\linewidth]{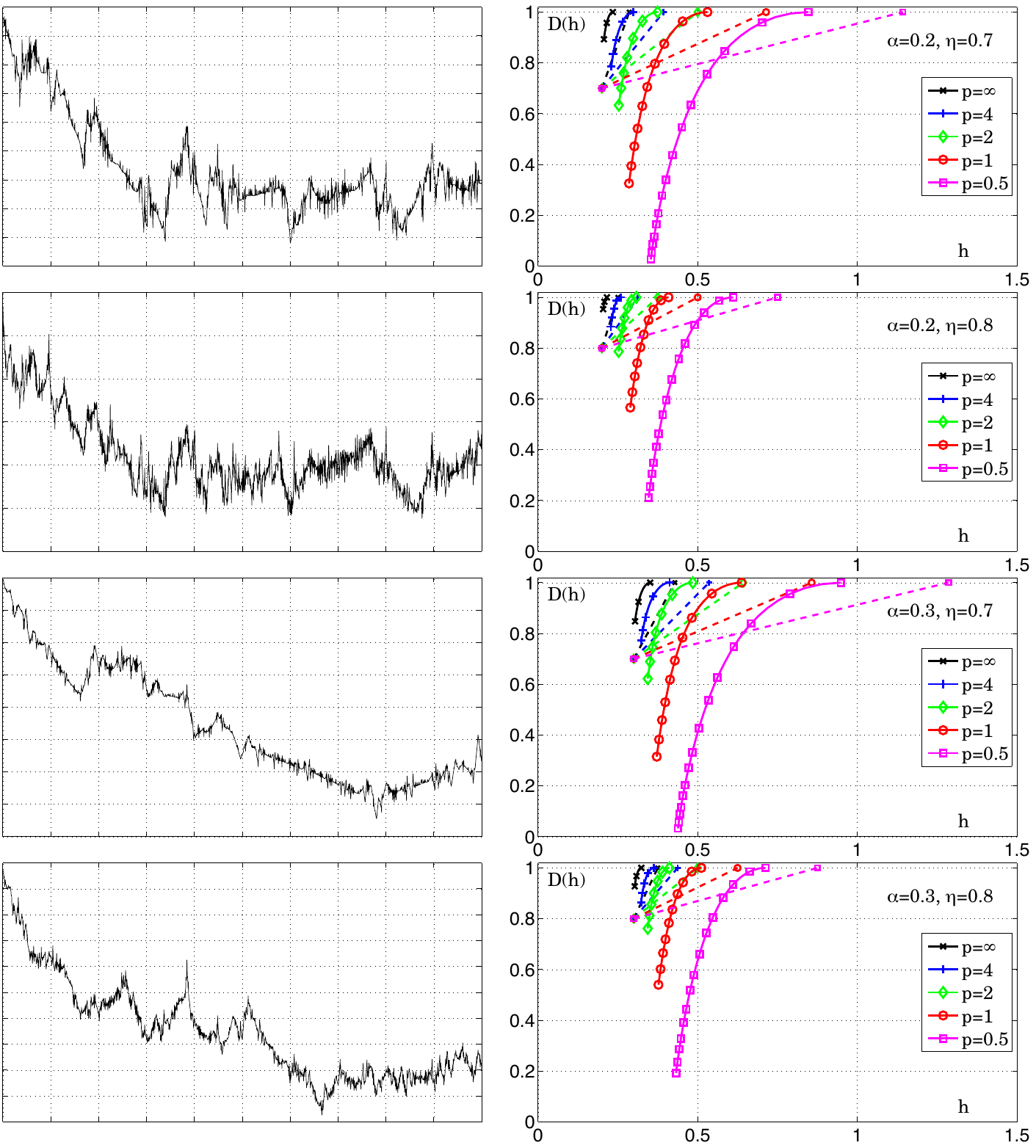}
\caption{\label{fig:estLWS}{\bf LWS} for $(\alpha,\eta)=\big\{(0.2,0.7),(0.2,0.8),(0.3,0.7),(0.3,0.8)\big\}$ (from top to bottom row): Single realizations (left column), theoretical spectra $\mathcal{D}(h)$ (right column, dashed lines) and estimated spectra $\mathcal{L}(h)$ (right column, solid lines).}
\end{figure}

\subsection{ $p$-spectra and lacunary singularities}
\label{sec:lacunary}

We illustrate the $p$-leader multifractal formalism for the estimation of the multifractal $p$-spectra of LWS, which contain lacunary singularities at almost every point and are hence not $p$-invariant.
Averages of $\mathcal{L}^{(p)}(h)$ over independent realizations are plotted in Fig. \ref{fig:estLWS} (colored solid lines with symbols), together with the theoretical spectra \eqref{eq:DhpLWS} (colored dashed lines) for four combinations of the parameters $(\alpha,\eta)$ with $\alpha\in\{0.2,0.3\}$ and $\eta\in\{0.7,0.8\}$.
Since the MFDFA method cannot reveal the difference in the $p$-spectra because it is limited to $p=2$, it is omitted in Fig. \ref{fig:estLWS}.

As expected (cf., \eqref{eq:DhpLWS}), the numerical estimates of the $p$-spectra are not invariant with $p$ but reproduce the evolution with $p$ of the theoretical spectra $\mathcal{D}^{(p)} (h)$: The larger $p$, the more the upper limit of the support of the spectra are shifted towards the point $\alpha$. The positions of the mode of $\mathcal{L}^{(p)}(h)$ slightly underestimate those of the true spectra, and the Hausdorff dimension of the leftmost point of the spectra are poorly estimated. Yet, the $\mathcal{L}^{(p)}(h)$ qualitatively reproduce the theoretical spectra satisfactorily well and, in particular, clearly and unambiguously reveal the lacunary nature of the sample paths.

\subsection{Images: Canonical Mandelbrot Cascades}
\label{sec:images}
As mentioned above, MFDFA has barely been used for images (except for the attempts in \cite{Gu2006,soares2009towards,wang2014local,shi2015new}), while the $p$-leader multifractal formalism extends in a straightforward manner to higher dimensions. Here, we illustrate this point and apply the $p$-leader multifractal formalism to synthetic multifractal images.

\noindent{\bf Canonical Mandelbrot Cascades (CMC).\quad}
The construction of multiplicative cascades of Mandelbrot (CMC) \cite{m74} is based on an iterative split-and-multiply procedure on an interval; we use a 2D binary cascade for two different multipliers: 
First, log-normal multipliers (CMC-LN), $W=2^{-U}$ with $U\sim \mathcal{N}(m,2m/\ln(2))$ a Gaussian random variable; Second, log-Poisson multipliers (CMC-LP), $W=2^\gamma\exp\left(\ln(\beta)\pi_\lambda\right)$, where $\pi_\lambda$ is a Poisson random variable with parameter $\lambda=-\frac{\gamma\ln(2)}{(\beta-1)}$.
We use fractional integration of order $\alpha=0.2$. CMC contain only canonical singularities, hence fractional integration results in a pure shift of their multifractal $p$-spectra by $\alpha$. Their multifractal $p$-spectra hence all collapse due to the $p$-invariance of canonical singularities. 
For CMC-LN, the  multifractal $p$-spectrum is given by
\[D^{(p)}(h)\equiv D(h)=2-\frac{(h-\alpha-m)^2}{4m}, 
\mbox{ with}\;  c_{1}=m+\alpha, c_2=-2m, c_m\equiv 0\; \mbox{ for all } \; m\geq3. \] 
The expression for the multifractal $p$-spectrum of CMC-LP is 
\[ {\cal D}(h)=2+\frac{\gamma}{\beta-1}+\frac{-\alpha+\gamma+h}{\ln \beta}\left[\ln\left(\frac{(-\alpha+\gamma+h)(\beta-1)}{\gamma\ln \beta}\right)-1\right], \] with 
$c_1=\alpha+\gamma\left(\frac{\ln(\beta)}{\beta-1}-1\right)$ and all higher-order log-cumulants are non-zero: 
$c_m=-\frac{\gamma}{\beta-1}\left(-\ln(\beta)\right)^m$. We set $m=0.04$, $\beta=0.8395$ and $\gamma=0.4195$, yielding $c_1=0.24$ and $c_2=-0.08$ for both cascades, and $c_3=0.014$ for CMC-LP. 

\begin{figure}[tb]
\centering%
\includegraphics[width=0.9\linewidth]{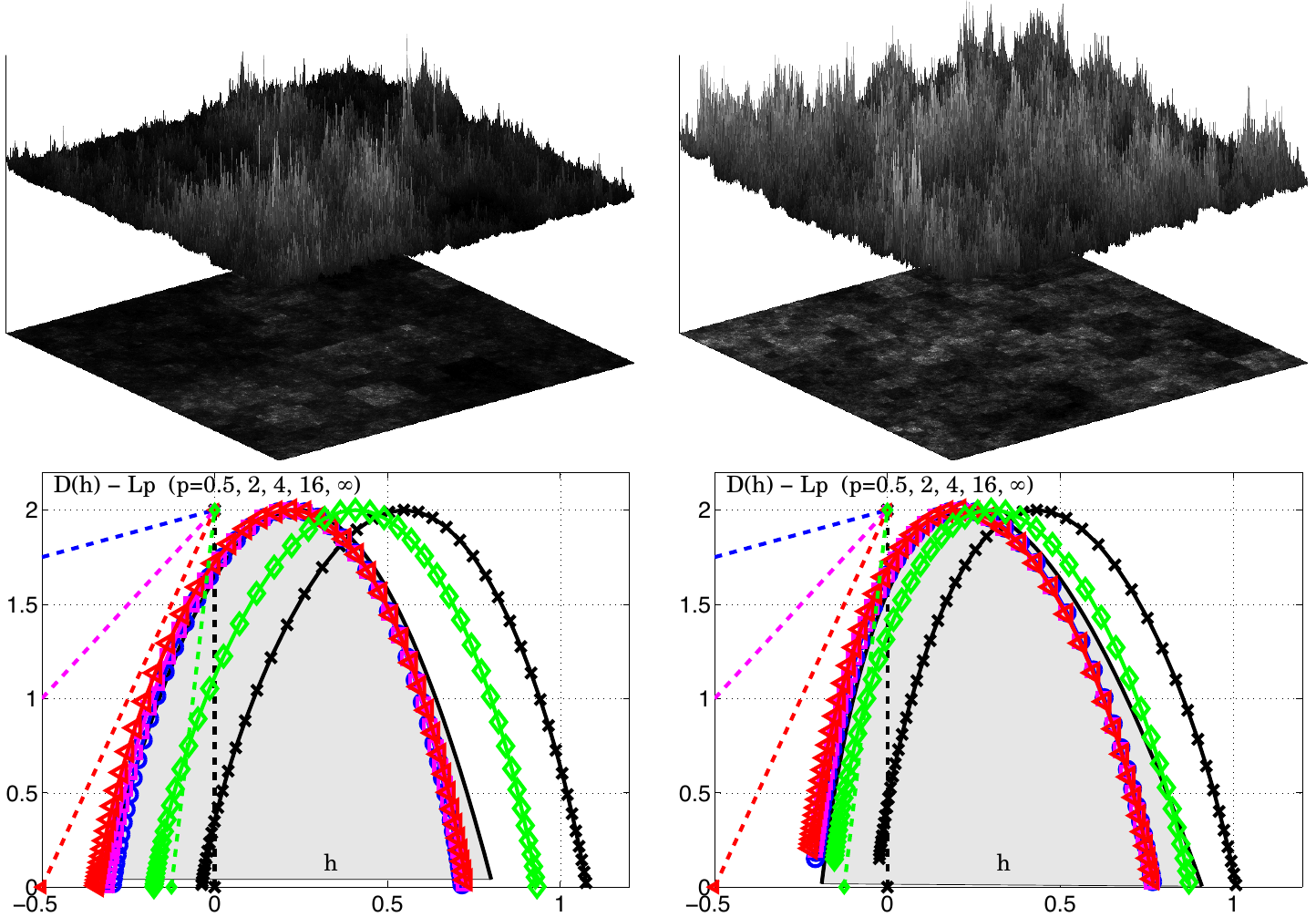}
\caption{\label{fig:CMC2D}{\bf 2D Mandelbrot cascade.} Log-Normal (left) and Log-Poisson (right) multipliers. Top: Single realization. Bottom: 
$\mathcal{D}(h)$ (black solid line and shaded area), $\mathcal{L}^{(p)}(h)$ for $p=+\infty$ (black, cross), $p=1/2$ (blue, circle), $p=2$ (magenta, square), $p=4$ (red, triangle), $p=16$ (green, diamond) and theoretical limits for multifractal $p$-spectra (dashed).}
\end{figure}

\noindent{\bf\boldmath Illustration of multifractal $p$-spectra.\quad} Averages over $100$ realizations of $p$-leader Legendre spectra $\mathcal{L}^{(p)}(h)$ 2D CMC of side length $N=2^{11}$ are reported in Fig. \ref{fig:CMC2D} (bottom) for CMC-LN (left) and CMC-LP (right), single realizations of the random fields are plotted in Fig. \ref{fig:CMC2D} (top); we use the scaling range $j_1=3$ and $j_2=8$ and tensor-product  Daubechies' wavelet with $N_\psi=2$ (cf. \cite{Antoine2004}).
The observations and conclusions are similar to those obtained in the previous subsection for (1D) signals (MRW); indeed, the $p$-leader estimates $\mathcal{L}^{(p)}(h)$ (solid lines in color, symbols) enable to correctly recover the theoretical $p$-spectrum $\mathcal{D}^{(p)}(h)$ as soon as the condition $p<p_0$ is fulfilled.

\section{Discussions, Conclusions and Perpectives}

The present contribution has developed and analyzed a novel multifractal formalism based on new local regularity exponents, the $p$-exponents, and their corresponding multiresolution quantities, the $p$-leaders, that have been theoretically defined and studied in the companion paper \cite{PART1}.
This new formulation of multifractal analysis generalizes the traditional H\"older-exponent-based formulation in several ways. First and foremost, it naturally allows to perform the multifractal analysis of functions that have negative regularity or, equivalently, that are not locally bounded (but instead belong locally to $L^p$). This allows to avoid the commonplace \emph{a priori} (fractional) integration of the data and the pitfalls it entails.

Moreover, the dependence on the parameter $p$ provides an important and rich information concerning the characteristics of the singular behavior of data. When the multifractal spectra differ for different valid values of $p$, this clearly indicates the presence of complex singular behaviors on the data, known as lacunary singularities. This information is a distinctive feature of $p$-exponents  and is not accessible with previous tools.

This contribution has also made a clear connection between the  $p$-leader multifractal formalism and MFDFA, a related technique that is widely used in applications. This has allowed to provide a theoretical grounding to MFDFA, which remained a useful but \emph{ad hoc} intuition. Also, it has brought to light that MFDFA actually measures the $2$-exponent, rather than the H\"older exponent as had been assumed previously. 

Numerical simulations on synthetic multifractal processes have shown that the $p$-leader  formalism, for small values of $p$, benefits from significant improvements in estimation performance, compared to wavelet leaders and MFDFA.
It is important to emphasize that, even though $p$-leaders were originally introduced with the goal of analyzing negative regularity, they show better estimation performance than wavelet leaders even for processes that have positive regularity only.  

To conclude, the benefits of the $p$-leader multifractal formalism over existing formulations are threefold: First, it allows the analysis of negative regularity; second, it shows better estimation performance; and third, the dependence of the estimates on the parameter $p$ provides richer and more detailed information on the characteristics of the singularities that can be  found in the data.
This suggests that the $p$-leader multifractal formalism, for small values of $p$, should be preferred for practical multifractal analysis.

A {\sc Matlab} toolbox implementing estimation procedures for performing $p$-multifractal analysis will be made publicly available at the time of publication.

\section*{Acknowledgements}

This work was supported by grants ANR AMATIS $\# 112432 $, 2012-2015, ANPCyT PICT-2012-2954 and PID-UNER-6136.

\appendix

  \section{Proofs}
  
  \label{proofs} 
  
  { \bf Proof of Proposition \ref{majpspec}}: We start by considering the case $p \geq 1$, in which case, we will prove the slightly sharper result that (\ref{bound:dh}) holds as soon as $f \in L^p$.   Let $H > -d/p$ be given. Let $E_H$ denote the set of points where the 
  $p$-exponent of $X$ is smaller than $H$. If $x_0 \in E_H$,  then there exists a sequence $r_n\rightarrow 0 $ such that 
  \[ \left( \frac{1}{r_n^d} \int_{B(x_0, r_n)} | f(x) |^p \; dx \right)^{1/p} \geq r_n^{H} ,\]
  so that there exists a sequence $j_n \rightarrow -\infty $ such that the dyadic cubes $\la_{j_n} (x_0)  $ satisfy 
  \[ \int_{ 3 \la_{j_n} (x_0)} | f(x) |^p \; dx \geq C 2^{j_n (d+Hp)} \]
(pick the smallest dyadic cube $\la_{j_n} (x_0)  $ such that $B(x_0, r_n)  \subset 3 \la_{j_n} (x_0)$). 
One of the $3^d$ cubes  of width $2^{j_n}$ that constitute $3 \la_{j_n} $ (which we denote by $ \mu_{j_n}$) satisfies 
 \[ \int_{ \mu_{j_n} (x_0)} | f(x) |^p \; dx \geq C3^{-d}  2^{j_n (d+Hp)} \]
 We consider now the maximal such dyadic cubes, of width less than a fixed  $\ep$,  satisfying this inequality for all possible $x_0$, and we denote by ${ \mathcal R}$ this collection. Then, since maximal dyadic  cubes  necessarily are 2 by 2 disjoint, 
 \[  C\sum_{\mu \in { \mathcal R} }   2^{j (d+Hp)} \leq \sum_{\mu \in { \mathcal R} }  \int_{ \mu} | f(x) |^p \; dx \leq C  \]
 where 
 $\mu$ is of width $2^j$. If $x\in E_H$, then $x$ belongs to one of the $3 \mu$, and therefore to the ball of same center and radius $3 d 2^{j}$. 
 Since $ r_n \geq C 2^{j_n}$, we  have obtained a covering of $E_H$ by balls of radius  at most $3 d \ep$ such that 
 \[ \sum diam (B(x, r))^{d+Hp}   \leq C,  \]
 and the result follows for $p \geq 1$. 
 
 We now consider the case the case $p < 1$. Hypothesis $\eta (p) >0$ means that 
 \BE \label{hypeta} \exists C,  \ep >0: \; \forall j <0, \quad 2^{dj} \sum_{i, k}  |c^{(i)}_{j,k}|^p \leq C 2^{\ep pj}. \EE 
 If $h_p (x_0) < H$, then there exists an infinite sequence of dyadic cubes $\la$ which contain $x_0$ and such that
\[   \sum_{j' \leq j, \;  \lambda' \subset 3 \lambda} \sum_{i=1}^{2^d-1} \big| c^{(i)}_{\la'} \big|^p \, \, 2^{-d(j-j')}  \geq 2^{H p j} . \] 
 We now consider the maximal cubes of width less that $\ep$ and which satisfy this condition. This yields a covering of $E_H$  by dyadic cubes $\la$ which are 2 by 2 disjoint. Denote by $A_j$ the cubes of this covering which are of  width $2^j$, and by $N_j$ the cardinality of $A_j$.  
 On one hand
 \[ \sum_\la \sum_{j' \leq j, \;  \lambda' \subset 3 \lambda} \sum_{i=1}^{2^d-1} \big| c^{(i)}_{\la'} \big|^p \, \, 2^{-d(j-j')} \geq N_j 2^{Hpj} \]
 (where the sum over $\la$ is taken on all dyadic cubes of width $2^j$). 
 On the other hand, the left side is equal to 
 \[ 3^d  \sum_{j' \leq j} \left( \sum_{ k' \in \ZZ^d} \sum_{i=1}^{2^d-1} \big| c^{(i)}_{\lambda'   } \big|^p \, \, 2^{-d(j-j')}\right).  \] 
 But (\ref{hypeta}) implies that the term between parentheses is bounded by $2^{-dj}2^{\ep p j'}$; thus we obtain that 
 \[ N_j 2^{Hpj} \leq C 2^{-dj} ,  \]  
 which implies that $\dim ( E_H) \leq d + Hp$,  and the result follows for $p < 1$. 

\bibliography{DFA2014}

\end{document}